\documentclass[9pt,twoside,english,reqno,a4paper]{amsart}
\usepackage{listings,graphicx,amsmath,varioref,amscd,amssymb,color,bm,stmaryrd,amsthm,amsfonts,graphics,geometry,latexsym,pgf,pst-all} 
\usepackage{amssymb, amsmath}
\usepackage{color}
\usepackage{palatino}
\definecolor{red}{rgb}{0.9,0,0}

\definecolor{azul}{rgb}{0,0,1}

\usepackage[colorlinks=false]{hyperref}
\usepackage{tikz-cd}
\usetikzlibrary{patterns}

\usepackage{indentfirst}
\usepackage{graphicx}

\theoremstyle{plain}
\newtheorem{theorem}{Theorem}[section]
\newtheorem{proposition}[theorem]{Proposition}

\newtheorem{lemma}[theorem]{Lemma}

\newtheorem{remark}[theorem]{Remark}
\newtheorem{definition}[theorem]{Definition}
\newtheorem{example}[theorem]{Example}

\def\re{\mathbb{R}}
\def\bk{\color{black}}
\def\dem{\noindent{\it Proof. }}

\newcommand{\fim}{\hfill $\Box$}

\usepackage{booktabs,amsmath}
\def\Xint#1{\mathchoice
    {\XXint\displaystyle\textstyle{#1}}%
    {\XXint\textstyle\scriptstyle{#1}}%
    {\XXint\scriptstyle\scriptscriptstyle{#1}}%
    {\XXint\scriptscriptstyle\scriptscriptstyle{#1}}%
      \!\int}
\def\XXint#1#2#3{{\setbox0=\hbox{$#1{#2#3}{\int}$}
    \vcenter{\hbox{$#2#3$}}\kern-.5\wd0}}
\def\dashint{\Xint-}

\title[Limit solutions to  $p-$Laplace problems involving Hardy potentials as $p\to1^+$]{Existence, non-existence and degeneracy of limit solutions to  $p-$Laplace problems involving Hardy potentials as $p\to1^+$}

\author[J. C. Ortiz Chata]{Juan Carlos Ortiz Chata}
\address{J. C. Ortiz Chata, Departamento de Matem\'atica, Universidade Federal de São Carlos - UFSCar, 13565-905, S\~ao Carlos - SP, Brazil, juanchata@ufscar.br}

\author[F. Petitta]{Francesco Petitta}
\address{Francesco Petitta, Dipartimento di Scienze di Base e Applicate per l' Ingegneria, Sapienza Universit\`a di Roma, Via Scarpa 16, 00161 Roma, Italia, francesco.petitta@uniroma1.it}

\date{}

\begin{document}
\pretolerance10000

\maketitle

\begin{abstract}
In this paper we analyze the asymptotic behaviour as $p\to 1^+$ of solutions $u_p$ to
$$
\left\{
\begin{array}{rclr}
-\Delta_p u_p&=&\frac{\lambda}{|x|^p}|u_p|^{p-2}u_p+f&\quad \mbox{ in } \Omega,\\
u_p&=&0 &\quad \mbox{ on }\partial\Omega,
\end{array}\right.
$$
where $\Omega$ is a bounded open subset of $\re^N$ with Lipschitz boundary, $\lambda\in\re^+$, and $f$ is a nonnegative datum in $L^{N,\infty}(\Omega)$. 
Under sharp smallness assumptions on the data $\lambda$ and $f$ we prove that $u_p$ converges to a suitable solution to the homogeneous Dirichlet problem
$$\left\{
\begin{array}{rclr}
					 \displaystyle - \Delta_{1} u &=& \frac{\lambda}{|x|}{\rm Sgn}(u)+f &  \text{in}\, \Omega, \\
					 u&=&0 & \text{on}\ \partial \Omega,			
				\end{array}\right.
$$
where $\Delta_{1} u ={\rm div}\left(\frac{D u}{|Du|}\right)$ is the $1$-Laplace operator. The main assumptions are  further discussed through explicit examples in order to show their optimality. 
\end{abstract}

\numberwithin{equation}{section}
\bibliographystyle{plain}
\maketitle

\section{Introduction}
Consider the problem
\begin{equation}
\label{Pintro} 
\left\{
\begin{array}{rclr}
-\Delta_1u&=&\frac{\lambda}{|x|}\frac{u}{|u|}+f&\quad \mbox{ in } \Omega,\\
u&=&0 &\quad \mbox{ on }\partial\Omega,
\end{array}
\right.
\end{equation} 
where $\Delta_1u=\hbox{div}\,\left(\frac{Du}{|Du|}\right)$, $\Omega\subset \mathbb{R}^N$ $(N\geq2)$ is an open set with bounded Lipschitz boundary containing the origin, $0<\lambda<N-1$, and $f$ belongs to $L^N(\Omega)$ satisfying suitable smallness assumptions (see \eqref{Hintro} below).

The solutions to this problem can be seen as the limit of the solutions of the following problems
\begin{equation}
\label{PLap}
\left\{
\begin{array}{rclr}
-\Delta_p u&=&\frac{\lambda}{|x|^p}|u|^{p-2}u+f&\quad \mbox{ in } \Omega,\\
u&=&0 &\quad \mbox{ on }\partial\Omega,
\end{array}
\right.
\end{equation}
as $p\to 1^+$.  The advantage of this auxiliary  problem is that it has been widely studied, for instance in  \cite{AzoreroAlonso1998}; here the authors studied critical elliptic and parabolic problems that include the study of \eqref{PLap}, which, in turn,  is obtained using variational approaches and providing  $\lambda$ is small enough and  $f$ belongs to the dual space $W^{-1,p}(\Omega)$. The case $p=2$ was then  studied from a purely PDE point of view in \cite{BoccardoOrsinaPeral2006} by showing existence and regularity of the solutions depending of the smallness of $\lambda$ and the summability of the datum $f$,  and it   was a first motivation for our study. 
We also mention that, in \cite{OrsinaOrtiz2024}, the authors obtained existence and optimal summability of solutions assuming that $f$ belongs to suitable Lebesgue's spaces $L^q(\Omega)$, which includes the case in which $q$ is close to 1. Several instances of problems involving Hardy type's terms as in \eqref{PLap} have been considered in the literature also in presence of lower order terms (see for instance  \cite{boc,lp, LMP2011, abo} and references therein).

On the other hand,  the $1-$Laplacian operator, $\Delta_1 u={\rm div}\,\left(\frac{Du}{|Du|}\right)$, is employed in variety of problems as the  image restoration model proposed by Rudin, Osher and Fatemi in \cite{ROF} in which the goal consists in the reconstruction  of an image $u$ given a blurry one $u_0=Au+n$, where A is a linear operator (e.g. a blur) and $n$ is a random noise.   This problem amount to look for minimizers of functionals as  
\begin{equation}\label{rof}
\int_{\Omega} |\nabla u| dx + \int_{\Omega}|A u- u_0|^2 dx,
\end{equation}
which are naturally relaxed in $BV(\Omega)$ in order to allow edges. Formally the sub-differential of the relaxed energy term in \eqref{rof} coincide with the 1-laplace operator.  

These type of operators  also enter, for instance,  in the study of torsional creep of  a cylindrical bar of constant cross section in $\re^2$, as well as in more theoretical issues of geometrical nature; for an incomplete account on these applications see for instance  \cite{OsSe, K, ka},  but also \cite{Sapiro, M, BCRS}, and  the monograph \cite{AndreuCasellesMazon2004}. 

The class of problems as  in \eqref{Pintro} when $\lambda=0$  have been studied in a series of papers \cite{K, ct, mst, MST2009} as an outcome of the study the behaviour of the solutions to the problem
\begin{equation}
\label{PLa}
\left\{
\begin{array}{rclr}
-\Delta_pu&=&f&\quad \mbox{ in } \Omega,\\
u&=&0 &\quad \mbox{ on }\partial\Omega,
\end{array}
\right.
\end{equation} 
as $p\to 1^+$ whenever the norm of $f$ is small. For instance, in \cite{K},  the author studied the existence of solutions for case $f=1$ provided $\Omega$ is suitably small. In particular, if \eqref{PLa} is formulated as a variational problem it is shown that it could admit  a non-trivial minimizer. In \cite{mst}, the authors dealt with the general case, when the datum belongs to $W^{-1,\infty}$, and in addition, they  proved that the limit of solutions for the above problem satisfies the concept of solution introduced in  \cite{acm2001, AndreuMazonMollCaselles2004, AndreuCasellesMazon2004}, also  showing that the smallness of the datum norm is essentially needed in order to get  the existence of a (non-trivial) solution. The case of problems with   $L^1$ data also in presence of lower order terms have been considered they also proved that the limit of its solutions satisfies this concept of solution (\cite{MST2009,  lops}).  

In \cite{Kawohl2007}, the authors studied Dirichlet problems for the 1-Laplace operator that including the eigenvalue problem and the following problem
$$
\left\{
\begin{array}{rclr}
-\Delta_1u&=&\lambda \frac{u}{|u|}+f&\quad \mbox{ in } \Omega,\\
u&=&0 &\quad \mbox{ on }\partial\Omega,
\end{array}
\right.
$$
where $f\in L^N(\Omega)$, for some $\lambda\in \mathbb{R}$. The approach to the problem above is variational, and it amount to look for minimizers of  a suitable relaxed  functional in $BV$. Moreover, with the aid of the duality theory and a general non-smooth Lagrange multiplier rule, it is found its Euler-Lagrange equation that contains a bounded vector field  ${\bf z}$, which can be identified with $Du/|Du|$ if $Du$ is nonzero and well-defined, and $u/|u|$, when $u(x)=0$, is substituted by some value in $[-1,1]$, which, in turn, corresponds  to the notion introduced in \cite{acm2001, AndreuMazonMollCaselles2004, AndreuCasellesMazon2004}. It is worth mentioning that among  the first works that deal with the eigenvalue problem for 1-Laplace operator  (say $f=0$) we have  \cite{Demengel2004, Demengel2005}. 

In \cite{ct}, the authors studied the problem \eqref{PLa} when the datum belongs to $L^N$ or to Marcikiewics space (Lorentz space) $L^{N,\infty}$, with small norm. In particular,  problem
\begin{equation}
\label{PC} 
\left\{
\begin{array}{rclr}
-\Delta_1u&=&\frac{\lambda}{|x|}&\quad \mbox{ in } \Omega,\\
u&=&0 &\quad \mbox{ on }\partial\Omega,
\end{array}
\right.
\end{equation} 
is covered. The fact that  $\frac{\lambda}{|x|}$ belongs to Lorentz space $L^{N,\infty}(\Omega)$, allowed  the authors to show the existence of non-trivial solution provided the datum, say $\lambda$, is suitably chosen. In \cite{dgop}, the authors studied the existence and regularity of non-trivial solutions for Dirichlet problems with a general singular nonlinearity, which include the problem \eqref{PC}.

The goal of our study is to provide a complete and optimal description of problem \eqref{Pintro} depending on the size of both the parameter $\lambda$ and the datum $f$; more precisely, if ${S_N}$ is the best constant of the Sobolev inequality for $W^{1,1}_0(\Omega)$, provided
\begin{equation}
\label{Hintro}
\|f\|_{L^N(\Omega)}S_N+\frac{\lambda}{N-1}\leq 1,
\end{equation}
we analyze the behaviour of the approximating solutions $u_p$ of problem \eqref{PLap} in order to show the existence of solutions for \eqref{Pintro}. 
 In order to do that, we need to introduce a suitable notion of solution for this problem and to make use of a  conveniently re-adapted Anzellotti-Chen-Frid theory of pairings between measure-divergence bounded vector fields and gradient of BV functions. If \eqref{Hintro} is satisfied with the strict inequality sign we show that $u_p\to 0$ as $p\to 1^{+}$ and $u=0$ is a solution of \eqref{Pintro} in the sense of our definition. Also be mean of explicit examples, in the extreme case in which 
$$
\|f\|_{L^{N,\infty}(\Omega)}\gamma+\frac{\lambda}{N-1}= 1, 
$$
($\gamma$ as in \eqref{ga}) we are able to show that solutions of problem \eqref{Pintro} obtained as limit of $u_p$ may be non-trivial, while non-existence of solution is shown to hold  provided 
$$
\|f\|_{L^N(\Omega)}S_N+\frac{\lambda}{N-1}> 1.
$$

This paper is organised as follows. In Section 2 we give some preliminaries tools, which  contain some rudiments on   Lorentz spaces,  functions of  bounded variations,  and an extended Anzellotti's theory. In Section 3 we present the main assumptions and results. In section 4 we study the existence of solutions to an approximate problem involving the $p-$Laplacian operator and give some basic estimative in order to pass to the limit as $p$ go to 1. In Sections 5 and 6, we prove the main theorems. In Section 7, we study the case for a more general data in $L^{N,\infty}(\Omega)$, or $W^{-1,\infty}(\Omega)$ and give some explicit examples.

\section{Preliminaries}

We fix the basic standard notation as follows: 
 \begin{itemize}
 \item $\mathcal{H}^{N-1}$ is the $N-1-$dimensional Hausdorff measure  on $\mathbb{R}^N$;
 \item $ \mathcal{L}^N$ is the $N-$dimensional Lebesgue measure on $\mathbb{R}^N$ and for each measurable set $E\subset \mathbb{R}^N$, we write $|E|$ for $\mathcal{L}^N(E)$;
 \item $B_r(x)$ is  the ball in $\mathbb{R}^N$ with center in $x$ and radius $r$;
  \item $\Omega$ is a open set in $\mathbb{R}^N$ and 
 $ \mathcal{D}(\Omega) (C^\infty_c(\Omega))$ is  the space of infinitely differentiable functions with compact support in $\Omega$;
  \item $\mathcal{D}'(\Omega)$ is  the space of distributions in $\Omega$;
 \item $ \mathcal{M}(\Omega,\mathbb{R}^N)$ is the space of the finite $\mathbb{R}^N-$valued Radon measure on $\Omega$.
 \end{itemize}
 
 Consider the set-valued sign function $\mbox{Sgn}\>:\>\mathbb{R}\to\mathcal{P}(\mathbb{R})$ given by
\begin{equation}
\label{Sgn}
\mbox{Sgn}(s)=
\left\{
\begin{array}{lcr}
1& \hbox{ if }& s>0\\
{[-1,1] }&\hbox{ if }& s=0\\
-1& \hbox{ if } &s<0
\end{array}
\right.
\end{equation}
and the truncation $T_k\>:\> \mathbb{R}\to \mathbb{R}$ by

$$
T_k(s)=
\left\{
\begin{array}{lcr}
s&\hbox{ if } &|s|\leq k\\
k\frac{s}{|s|}&\hbox{ if } &|s|>k.
\end{array}
\right.
$$
\medskip
\subsection{Some inequalities in Sobolev spaces}

Let $1\le p< N$ and $p^*=Np/(N-p)$. 

\begin{proposition}[Sobolev's inequality]
There exists a constant $C=C(N,p)>0$ such that, for any $ u\in W^{1,p}(\mathbb{R}^N) $, 
\begin{equation}
\label{SoI}
\left(\int_{\mathbb{R}^N}|u|^{p^*}\,dx\right)^\frac{1}{p^*}\leq C\left(\int_{\mathbb{R}^N}|\nabla u|^p\,dx\right)^\frac{1}{p}.
\end{equation}
\end{proposition}
We denote by $S_{N, p}$ the better constant in \eqref{SoI}, i.e.,
$$
S_{N,p}=\sup_{u \in W^{1,p}(\mathbb{R}^N)\setminus\{0\}}\frac{\|u\|_{L^{p^*}(\mathbb{R}^N)}}{\|u\|_{W^{1,p}(\mathbb{R}^N)}}.
$$
For $p=1$, we write $S_{N}$ for $S_{N,1}$. More details about $S_{N, p}$, it can be found in \cite{Talenti}.

\begin{proposition}[Hardy's inequality]\label{Halp} For any $u\in W^{1,p}(\mathbb{R}^N)$ the following inequality is established
\begin{equation}
\label{HaI}
\int_{\mathbb{R}^N}\frac{|u|^p}{|x|^p}\,dx\leq \frac{1}{\mathcal{H}^p}\int_{\mathbb{R}^N}|\nabla u|^p\,dx,
\end{equation}
where $\mathcal{H}=\frac{N-p}{p}$. Furthermore, $\mathcal{H}$ is the optimal constant in \eqref{HaI}.
\end{proposition}
 Proposition \ref{Halp} is a particular instance  of the so-called Caffarelli-Kohn-Nirenberg inequality \cite{CKN1984} which gives that  there exists a constant $C>0$ such that
$$
\left(\int_{\mathbb{R}^N}\frac{1}{|x|^p}|\varphi|^p\,dx\right)^\frac{1}{p}\leq C\left(\int_{\mathbb{R}^N}|\nabla \varphi|^p\,dx\right)^\frac{1}{p}\quad\hbox{ for all }\; \varphi\in C^\infty_c(\mathbb{R}^N).
$$
Since $C^\infty_c(\mathbb{R}^N)$ is dense in $W^{1,p}_0(\mathbb{R}^N)=W^{1,p}(\mathbb{R}^N)$, this inequality holds for all $u\in W^{1,p}(\mathbb{R}^N)$.

In fact  that $C=\left(\frac{p}{N-p}\right)^p$ is the optimal constant for all $p\in [1, N)$, recalling that for $p\in (1, N)$, is  proven in  \cite[Lemma 2.1]{AzoreroAlonso1998}. For $p=1$, it follows, by an argument  given in \cite{Talenti}, that
$$
C=\frac{1}{N-1}\,.
$$

\subsection{ Lorentz spaces}

Let $f$ be a measurable function in $\Omega$. We recall that its distribution function is given by 
$$
\alpha_f(s)=|\{x\in \Omega\>: \> |f(x)|>s\}|
$$
and that its non-increasing rearragement in $(0, +\infty)$ is defined as
$$
\label{NRf}
f^*(t)=\inf\{s>0\>:\> \alpha_f(s)\leq t\}.
$$
As usual we let $\inf \emptyset=+\infty$.

\begin{definition}
For $1\leq p <\infty$ and $1\leq q \leq \infty$, the Lorentz space $L^{p, q}(\Omega)$ is defined as 
$$
L^{p,q}(\Omega)=\{ f\>:\> f\mbox{ measurable on }\>\Omega,\; \|f\|_{L^{p,q}(\Omega)}<\infty\},
$$
where 
$$
\|f\|_{L^{p,q}(\Omega)}=\left\{
\begin{array}{lcr}
\left[\int_{0}^{\infty}[t^\frac{1}{p}f^*(t)]^q\frac{dt}{t}\right]^\frac{1}{q}& \hbox{ if } &1\leq p<\infty,\; 1\leq q<\infty\\
\sup_{t>0}t^\frac{1}{p}f^*(t)& \hbox{ if } & 1\leq p< \infty,\; q=\infty.
\end{array}
\right.
$$

\end{definition}

\begin{proposition}[H\"older's inequality]
\label{HI}
Let $1<p<\infty$, $ 1\leq q \leq  \infty$, $1/p+1/p'=1$ and $1/q+1/q'=1$. Let $f\in L^{p,q}(\Omega)$ and $g\in L^{p',q'}(\Omega)$. Then 
$$
\left|\int_{\Omega}f(x)g(x)\,dx\right|\leq \|f\|_{L^{p, q}(\Omega)}\|g\|_{L^{p', q'}(\Omega)}.
$$
\end{proposition}
The proof of the previous  can be found in \cite{ONeil1968} (see also \cite[Appendix]{DAnconaFanelli2007}).

\subsection{The space $BV$}

Recall that the space of functions of bounded variation in $\Omega$, $BV(\Omega)$, is made up of all functions of $L^1(\Omega)$ whose distributional derivative is a finite Radon measure in $\Omega$, that is,
$$
BV(\Omega)=\{u\in L^1(\Omega)\>;\> Du\in \mathcal{M}(\Omega, \mathbb{R}^N)\}.
$$
Furthermore $u$ is a bounded variation function in $\Omega$ if and only if $u\in L^1(\Omega)$ and  the following quantity is finite 
\begin{equation}
\label{TV}
\sup\left\{\int_{\Omega}u\,\mbox{div}\,(\varphi)\,dx\>:\>\varphi\in C^1_c(\Omega,\mathbb{R}^N),\quad\|\varphi\|_{L^\infty(\Omega)}\leq 1\right\},
\end{equation}
 which, in turn, by \cite[Proposition 3.6]{AmbrosioFuscoPallara}, is equal to the total variation of $Du$ on $\Omega$, $\displaystyle |Du|(\Omega)=\int_{\Omega}|Du|$.

The space $BV(\Omega)$ is a Banach space endowed with the norm
$$
\|u\|_{BV(\Omega)}=\int_{\Omega}|Du|+\int_{\Omega}|u|\,dx.
$$

On the other hand, if $\Omega\subset \mathbb{R}^N$ is an open set with bounded Lipschitz boundary, the classical result of existence of trace operator $u\mapsto u^\Omega$ between Sobolev spaces $W^{m,p}(\Omega)$ and $L^1(\partial\Omega)$ can be extended between $BV(\Omega)$ and $L^1(\partial\Omega)$ (for instance see \cite{AmbrosioFuscoPallara, AttouchButtazzoMichaille2006, EvansGariepy}). In particular, it is established the following:

\begin{proposition}
\label{BT}
Let $u\in BV(\Omega)$. Then, for $\mathcal{H}^{N-1}-$almost every $x\in \partial\Omega$ there exists $u^\Omega(x)\in \mathbb{R}$ such that
$$
\lim_{\rho\to 0}\rho^{-N}\int_{\Omega\cap B_\rho(x)}|u(y)-u^\Omega(x)|\,dy=0.
$$

Moreover, 
$$
\|u^\Omega\|_{L^1(\partial\Omega)}\leq C\|u\|_{BV(\Omega)},
$$
for some constant $C$ depending only on $\Omega$.
\end{proposition}

\begin{definition}
Let $u, u_n\in BV(\Omega)$, for all $n\in \mathbb{N}$, we say that $(u_n)_{n\in\mathbb{N}}$ strictly converges in $BV(\Omega)$ to $u$ if
\begin{eqnarray*}
u_n\to u\quad \hbox{ in }\; L^1(\Omega),\\
\int_{\Omega}|Du_n|\to \int_{\Omega}|Du|,
\end{eqnarray*}
as $n\to \infty$.
\end{definition}

\begin{proposition}
Let $\Omega\subset \mathbb{R}^N$ be an open set with bounded Lipschitz boundary. Then, the trace operator $u\mapsto u^\Omega$ is continuous between $BV(\Omega)$ and $L^1(\partial\Omega)$ with the topology induced by strict convergence.
\end{proposition}
For simplicity from now on we will write $u$ ($\in L^1(\partial\Omega)$) instead of $u^\Omega$.

\begin{remark}
Since the trace operator $ u\mapsto u\in L^1(\partial\Omega)$ is continuous and onto and  the embedding $BV(\Omega)\hookrightarrow L^{1^*}(\Omega)$ is continuous,  the norm above is equivalent to the following
$$
\|u\|:=\int_{\Omega}|Du|+\int_{\partial\Omega}|u|\,d\mathcal{H}^{N-1}.
$$
 \end{remark}

\begin{proposition}
Let $(u_n)_{n\in\mathbb{N}}$ be a sequence in $BV(\Omega)$ such that $\sup_{n\in \mathbb{N}}|Du_n|(\Omega)<\infty$ and $u_n\to u$ in $L^1(\Omega)$. Then one has that $u\in BV(\Omega)$ and 
$$
\int_{\Omega}|Du|\leq\liminf_{n\to\infty}\int_{\Omega}|Du_n|.
$$
\end{proposition}
\dem
It follows immediately by observation in \eqref{TV}.
\fim 

Moreover the following inequality is stablished in $BV(\Omega)$.
\begin{proposition}
\label{PI}
Let $\Omega\subset \mathbb{R}^N$ be an open set with bounded Lipschitz boundary. Let $0\leq b\leq1$. Then there exists a constant $C>0$ such that the following inequality holds
$$
\left(\int_{\Omega}\frac{1}{|x|^{b1^*_b}}|u|^{1^*_b}\,dx\right)^\frac{1}{1^*_b}\leq C\left(\int_{\Omega}|Du|+\int_{\partial\Omega}|u|\,d\mathcal{H}^{N-1}\right)\quad\hbox{ for all }\;u\in BV(\Omega),
$$
where $1^*_b=\frac{N}{N-(1-b)}$.
\end{proposition}

The proof for the case $b=0$ is given in \cite[Theorem 3.47]{AmbrosioFuscoPallara} and for $0<b<1$ in \cite{OrtizPimentaSegura2021}.

\begin{proposition}[Sobolev's inequality]
\label{LSI}
There exists a constant $\gamma>0$ such that
\begin{equation}\label{bestl}
\|u\|_{L^{1^*, 1}(\mathbb{R}^N)}\leq \gamma \int_{\mathbb{R}^N}|Du| \quad\mbox{ for all }\; u\in BV(\mathbb{R}^N).
\end{equation}
The better constant is given by
$$
\gamma = \frac{1}{(N-1)|B_1(0)|^\frac{1}{N}}.
$$
\end{proposition} 

In \cite{Alvino1977}, the author showed that if $u$ is a sufficiently smooth function then the following inequality holds
$$
\|u\|_{L^{1^*, 1}(\mathbb{R}^N)}\leq \gamma\int_{\mathbb{R}^N}|Du|.
$$
Hence, using an approximation argument as in \cite{EvansGariepy, Ziemer1989}, one may prove that this inequality is true for all $u\in BV(\mathbb{R}^N)$.

\subsection{ Some fine properties of $BV$ functions}

Recall that, for $u\in L^1(\Omega)$, $u$ has an approximate limit at $x\in \Omega$ if there exists $\tilde{u}(x)$ such that 
$$
\lim_{\rho\downarrow0} \dashint_{B_{\rho}(x)}|u(y)-\tilde{u}(x))|\,dx=0,
$$
where $\dashint_{E}f\,dx=\frac{1}{|E|}\int_{E}f\,dx$. The set where this property does not hold is denoted by $S_u$. This is a $\mathcal{L}^N-$negligible Borel set \cite[Proposition 3.64]{AmbrosioFuscoPallara} . We say that $x$ is an approximate jump point of $u$ if there exists $u^+(x)\neq u^{-}(x)$ and $ \nu\in S^{N-1}$  such that
\begin{eqnarray*}
\lim_{\rho \downarrow 0}\dashint_{B^+_\rho(x, \nu) }|u(y)-u^+(x)|\,dx&=&0\\
\lim_{\rho \downarrow 0}\dashint_{B^{-}_\rho(x, \nu) }|u(y)-u^{-}(x)|\,dx&=&0,
\end{eqnarray*}
where
\begin{eqnarray*}
B^+_\rho (x, \nu)&=&\{y\in B_\rho(x)\>:\> \langle y-x, \nu\rangle >0\}\\
B^{-} _\rho (x, \nu)&=&\{y\in B_\rho(x)\>:\> \langle y-x, \nu\rangle <0\}.
\end{eqnarray*}
The set of approximate jump points is denoted by $J_u$. The set $J_u$ is a Borel subset of $S_u$ \cite[Proposition 3.69]{AmbrosioFuscoPallara} and $\mathcal{H}^{N-1}(S_u \setminus J_u) = 0$.

For $u\in L^1(\Omega)$, $u^*\>:\>\Omega\setminus(S_u\setminus J_u)\to \mathbb{R}$ is called the precise representative of u if 
\begin{equation}
\label{CanRep}
u^*(x)=
\left\{
\begin{array}{lcr}
\tilde{u}(x)&\mbox{ if }& x\in \Omega\setminus S_u\\
\frac{u^+(x)+u^{-}(x)}{2}& \mbox{ if }& x\in J_u.
\end{array}
\right.
\end{equation}

\begin{proposition}
Let $u \in BV(\Omega)$. The mollified functions $u*\rho_\epsilon$ pointwise converge to $u^*$ in its domain.
\end{proposition}

The proof of the above proposition as well as to more details about space $BV(\Omega)$ can be found in \cite{AmbrosioFuscoPallara} (to which we refer, in general, for further standard  notations not explicitly recalled here), see also \cite{AttouchButtazzoMichaille2006, EvansGariepy}.

\subsection{An extension of Anzellotti's theory}

In this subsection, we recall an extension of Anzellotti's theory \cite{Anzellotti1983}, which is mainly outlined in  \cite{Caselles2011}.

Let $\Omega\subset \mathbb{R}^N$ be an open set with bounded Lipschitz boundary. We set the space
$$
X_{\mathcal{M}}(\Omega)=\{{\bf z}\in L^{\infty}(\Omega,\mathbb{R}^N)\>:\> \mbox{div}\,{\bf z} \; \mbox{ is a finite Radon measure in }\;\Omega\}.
$$

We define
$$  
\langle \mbox{div}\,{\bf z},u\rangle=\int_{\Omega}u^*\,\mbox{div}\,{\bf z}\quad\hbox{ for all }\; {\bf z}\in X_{\mathcal{M}}(\Omega),\;\hbox{ for all }\; u\in BV(\Omega)\cap L^\infty(\Omega),
$$
where $u^*$ is as in \eqref{CanRep}.

For each $u\in W^{1,1}(\Omega)\cap L^\infty(\Omega)$ and ${\bf z}\in X_{\mathcal{M}}(\Omega)$ we define
$$
\langle {\bf z}, u\rangle_{\partial\Omega}=\int_{\Omega}u^*\,\mbox{div}{\bf z}+\int_{\Omega}{\bf z}\cdot \nabla u\,dx.
$$
\begin{proposition}
If $u, v \in W^{1,1}\cap L^\infty(\Omega)$ and $ u = v $ on $ \partial\Omega $ then
$$
\langle {\bf z}, u\rangle_{\partial\Omega} = \langle {\bf z}, v\rangle_{\partial\Omega}.
$$
\end{proposition}
\dem Since $ u-v\in W^{1,1}_0(\Omega) $, there exist $ (\varphi_m) \subset C^\infty_c(\Omega) $ such that
$$
\varphi_m\to u-v\quad\hbox{ in }\; W^{1,1}(\Omega)\quad\hbox{ as }\; m\to \infty.
$$
Thus
\begin{align*}
\langle {\bf z}, u-v \rangle_{\partial\Omega}&= \int_{\Omega}(u-v)^*\hbox{ div }{\bf z}+\int_{\Omega}{\bf z}\cdot \nabla (u-v)\,dx\\
&=\lim_{m\to\infty}\left(\int_{\Omega}\varphi_m\hbox{div}{\bf z}+\int_{\Omega}{ \bf z}\cdot \nabla \varphi_m\,dx\right)\\
&=0.
\end{align*} 
\fim
\begin{definition}
For each ${\bf z}\in X_{\mathcal{M}}(\Omega)$ and $u\in BV(\Omega)$, we define
$$
\langle {\bf z}, u\rangle_{\partial\Omega}=\langle {\bf z}, w\rangle_{\partial\Omega},
$$
where $w\in W^{1,1}(\Omega)\cap L^\infty(\Omega)$ is such that $ w = u $ on $ \partial\Omega $. 
\end{definition}

The definition above is well defined since the trace operator from $BV(\Omega)$ in $L^1(\partial\Omega)$ is onto and \cite[Lemma 5.5]{Anzellotti1983} holds. On the other hand, arguing as in the proof of \cite[Theorem 1.2]{Anzellotti1983} we get the following resulted.

\begin{proposition}
\label{Prop2}
Let ${\bf z}\in X_{\mathcal{M}}(\Omega)$. Then there exists a linear operator $ \gamma\>:\>X_{\mathcal{M}}(\Omega)\to L^\infty(\partial\Omega)$ such that 
\begin{eqnarray*}
\label{E1}
\|\gamma_{\bf z}\|_{ L^\infty(\partial\Omega)}&\leq&\|{\bf z}\|_{L^\infty(\Omega,\mathbb{R}^N)}\\
\langle {\bf z}, u\rangle_{\partial\Omega}&=&\int_{\partial\Omega}\gamma_{\bf z}(x)u(x)\,d\mathcal{H}^{N-1}\quad \mbox{ for all }u\in BV(\Omega) \cap L^\infty(\Omega)\\
\label{E3}
\gamma_{\bf z}(x)&=&{\bf z}(x)\cdot \nu(x)\quad\mbox{ for all }\; x\in \partial\Omega \quad\mbox{ if }\; {\bf z}\in C^1(\overline{\Omega},\mathbb{R}^N).
\end{eqnarray*}
\end{proposition}

The function $\gamma_{\bf z}(x)$ is a weakly defined trace on $\partial\Omega$ of the normal component of ${\bf z}$, which will be denoted by $[{\bf z}, \nu]$. 

\begin{definition}
Let ${\bf z} \in X_{\mathcal{M}}(A)$ and $u\in BV(A)\cap L^\infty(A)$ for all open set $A\subset\subset\Omega$. We define a linear functional $({\bf z}, Du)\>:\>C^\infty_c(\Omega)\to \mathbb{R}$ as 
$$
\langle ({\bf z}, Du),\varphi\rangle=-\int_{\Omega}u^*\varphi\,\mbox{div}\,{\bf z}-\int_{\Omega}u {\bf z}\cdot \nabla \varphi\,dx.
$$
\end{definition}

\begin{proposition}
\label{Dist}
Let $A\subset \Omega$ be a open set. For all functions $\varphi\in C_c(A)$ the following inequality holds                                                                                                                                                                                                                                                                                                                                                                                                                
$$
|\langle ({\bf z}, Du), \varphi\rangle |\leq \sup_{A} |\varphi|\|{\bf z}\|_{L^\infty(A, \mathbb{R}^N)}\int_{A}|Du|,
$$
which means that $({\bf z}, Du)$ is a Radon measure in $\Omega$. Furthermore, $({\bf z}, Du)$, $|({\bf z}, Du)|$ are absolutely continuous with respect to the measure $|Du|$ in $\Omega$ and 
$$
\left|\int_{B}({\bf z}, Du)\right|\leq \int_{B} |({\bf z}, Du)|\leq \|{\bf z}\|_{L^\infty(A,\mathbb{R}^N)}\int_{B}|Du|\,,
$$
for all Borel sets $B$ and for all opens sets $A$ such that $B\subset A\subset \Omega$.
\end{proposition}

\begin{proposition}
\label{Mea}
Let ${\bf z}\in X_{\mathcal{M}}(\Omega)$ and $u\in BV(\Omega)\cap L^\infty(\Omega)$. Then the following identity holds
$$
\int_{\Omega}u^*\, \mbox{div}\,{\bf z}+\int_{\Omega}({\bf z}, Du)=\int_{\partial\Omega}[{\bf z}, \nu]\,u\,d\mathcal{H}^{N-1}.
$$
\end{proposition}
Proofs of the propositions \ref{Dist} and \ref{Mea} can be found in \cite[Section 5]{Caselles2011} (see also \cite{Anzellotti1983}).

\begin{definition}
Let ${\bf z}\in X_{\mathcal{M}}(\Omega)$ and $u\in BV(\Omega)$ be such that $u^*\in L^1(\Omega, |\hbox{div}\, {\bf z}|)$. The linear functional $({\bf z}, Du)\>:\>C^\infty_c(\Omega)\to \mathbb{R}$ is defined by
$$
\langle ({\bf z}, Du), \varphi\rangle=-\int_{\Omega}\varphi u^*\hbox{div}\,{\bf z}-\int_{\Omega} u {\bf z}\cdot \nabla \varphi\,dx\quad\hbox{ for all }\, \varphi\in C^\infty_c(\Omega).
$$
\end{definition}

\begin{proposition}
\label{Tku}
Suppose that ${\bf z}\in X_{\mathcal{M}}(\Omega)$ and $u\in BV(\Omega)$ is such that $u^*\in L^1(\Omega, |\hbox{div}\,{\bf z}|)$ and $u^+, u^-\in L^1(J_u, |\hbox{div}\,{\bf z}|)$. Then
$$
({\bf z}, T_k(u))\to ({\bf z}, Du)\quad\hbox{ in }\; \mathcal{D}'(\Omega)\quad\hbox{ as }\; k\to \infty.
$$
\end{proposition}
\dem
In order to proof this convergence, it is sufficient to prove that
$$
\int_{\Omega}\varphi T_k(u)^*\hbox{div}\,{\bf z}\to \int_{\Omega}\varphi u^*\hbox{div}\,{\bf z}\quad\hbox{ as}\;k\to \infty.
$$
By \cite[Proposition 3.64, 3.69]{AmbrosioFuscoPallara}, we can write
\begin{align}
\label{EstTk}
\nonumber\int_{\Omega}\varphi T_k(u)^*\hbox{div}\,{\bf z}&=\int_{\Omega\setminus S_u}\varphi T_k(\tilde u)\hbox{div}\,{\bf z}+\int_{J_u\cap J_{T_k(u)}}\frac{T_k(u^+)+T_k(u^-)}{2}\varphi\hbox{div}\,{\bf z}\\ 
&\qquad+\int_{J_u\setminus  J_{T_k(u)}}  k {\rm sgn}(u^+) \varphi\hbox{div}\,{\bf z}.
\end{align}
{\bf Claim 1:}
\begin{equation*}
\lim_{k\to\infty}\int_{J_u}\frac{T_k(u^+)+T_k(u^-)}{2}\chi_{J_{T_k(u)}}(x)\hbox{div}\,{\bf z}=\int_{J_u}\frac{u^++u^-}{2}\hbox{div}\,{\bf z}.
\end{equation*}
Indeed, let $x_0\in J_u$, there exists $k_0\geq 1$ such that
$$
|u^\pm(x_0)|\leq k_0.
$$
Thus, for each $k\geq k_0$, we have $T_k(u^+(x_0))\neq T_k(u^-(x_0))$ and so
$$
x_0\in J_{T_k(u)}\quad\hbox{ for all }\; k\geq k_0.
$$
Hence, 
$$
\frac{T_k(u^+(x_0))+T_k(u^-(x_0))}{2}\chi_{J_{T_k(u)}}(x_0)=\frac{u^+(x_0)+u^-(x_0)}{2}\quad\hbox{ for all }\; k\geq k_0.
$$ 
and, since $u^+, u^-\in L^1(J_u, |\hbox{div}\,{\bf z}|)$, by Lebesgue's dominated convergence theorem, Claim 1 holds.

{\bf Claim 2: }
\begin{equation*}
\int_{J_u\setminus J_{T_k(u)}} k {\rm sgn}(u^+)\varphi\hbox{div}\,{\bf z}\to 0\quad\hbox{ as }\; k\to \infty.
\end{equation*}
We write
$$
\int_{J_u\setminus J_{T_k(u)}} k {\rm sgn}(u^+)\varphi\hbox{div}\,{\bf z}=\int_{J_u\setminus S_{T_k(u)}} k {\rm sgn}(u^+)\varphi\hbox{div}\,{\bf z}+\int_{(S_{T_k(u)}\setminus J_{T_k(u)})\cap J_u}k {\rm sgn}(u^+)\varphi\hbox{div}\,{\bf z}.
$$
Since $\mathcal{H}^{N-1}(S_{T_k(u)}\setminus J_{T_k(u)})=0$, by \cite[Proposition 3.1]{ChenFried1999}, the second term on the right-hand side above vanishes. 

In order to estimate the first term on the right-hand side above, note that, for $l\geq k\geq 1$, by \cite[Proposition 3.64]{AmbrosioFuscoPallara}, $S_{T_k(u)}\subset S_{T_l(u)}$. Thus
$$
\lim_{k\to \infty}\int_{J_u\setminus S_{T_k(u)}}\max\{|u^+|, |u^-|\}|\varphi||{\rm div}\,{\bf z}|=\int_{J_u\setminus  \bigcup_{k\geq 1}S_{T_k(u)}}\max\{|u^+|, |u^-|\}|\varphi||{\rm div}\,{\bf z}|=0,
$$
since $\bigcup_{k\geq1}S_{T_k(u)}\setminus (S_u\setminus J_u)=J_u$ and $u^+, u^-\in L^1(J_u, |\hbox{div}\,{\bf z}|)$. But this, in turn, implies, using the fact that $k< \max\{ |u^+|, |u^-|\}$ on $J_u\setminus J_{T_k(u)}$, that
$$
\limsup_{k\to\infty}\left|\int_{J_u\setminus S_{T_k(u)}}k {\rm sgn}\,(u^+)\varphi{\rm div}\,{\bf z}\right|\leq \lim_{k\to\infty}\int_{J_u\setminus S_{T_k(u)}}\max\{|u^+|, |u^-|\}|\varphi||{\rm div}\,{\bf z}|=0.
$$

Hence Claim 2 holds.

Therefore, passing to the limit in  \eqref{EstTk}, by Claim 1 and Claim 2, we get that
$$
\lim_{k\to \infty}\int_{\Omega}\varphi(T_k(u))^*\hbox{div}\,{\bf z}=\int_{\Omega\setminus J_u}\varphi \tilde{u}\,\hbox{div}\,{\bf z}+\int_{J_u}\varphi\frac{u^++u^-}{2}\hbox{div}\,{\bf z}=\int_{\Omega}\varphi u^*\,\hbox{div}\,{\bf z}\,,
$$
as desired.
\fim

\section{Main assumptions and results}\label{main}

The main focus  of this paper is on the existence and nonexistence  of solutions to the problem
\begin{equation}
\label{P} 
\left\{
\begin{array}{rclr}
-\Delta_1u&=&\frac{\lambda}{|x|} \frac{u}{|u|}+f&\quad \mbox{ in } \Omega,\\
u&=&0 &\quad \mbox{ on }\partial\Omega,
\end{array}
\right.
\end{equation} 
where  $ \Delta_1u=\mbox{div}\,\left(\frac{Du}{|Du|}\right)$ and $\Omega\subset \mathbb{R}^N$ be an open set with bounded Lipschitz boundary containing the origin, $ \lambda\in(0,N-1)$ and, as a starting point,  $f\in L^N(\Omega)$ is a not identically null function. We will mainly discuss, as a structural assumption, the following: 
\begin{equation}
\label{H}
\|f\|_{L^N(\Omega)}S_N+\frac{\lambda}{N-1}\leq 1.
\end{equation}

As already seen, the above problem is a borderline case with respect of previous works involving the $1-$Laplacian operator \cite{K, ct, mst} due to the presence of the zero order perturbation term 
$$
\frac{\lambda}{|x|}\frac{u}{|u|},
$$
as $1/|x|$ belongs to Marcinkiewicz (or Lorentz) space $L^{N,\infty}(\Omega)$ but not to $L^N(\Omega)$. The main feature of this work is to analyse, assuming \eqref{H}, the influence of this perturbation term on the existence of a solutions for \eqref{P}. Moreover, as $Du/|Du|$ does not mean when $Du=0$, a precise definition of solution for \eqref{P} is needed, which will be given in sense of the concept of solution introduced in \cite{acm2001, AndreuMazonMollCaselles2004, AndreuCasellesMazon2004}. 

 \begin{definition}
\label{DefSol}
We say that $u\in BV(\Omega)$ is a solution of \eqref{P} if  there exist ${\bf z}\in X_{\mathcal{M}}(\Omega)$ with $\|{\bf z}\|_{L^\infty(\Omega,\mathbb{R}^N)}\leq 1$ and $s(x)\in \hbox{Sgn}(u(x))$ a.e. $x\in\Omega$ such that
\begin{itemize}
\item[$(1)$] $-\mbox{div}\,{\bf z}=\frac{\lambda}{|x|}s(x)+f\quad \mbox{ in }\;\mathcal{D}'(\Omega)$;
\item[$(2)$] $({\bf z}, DT_k(u))=|DT_k(u)|\quad\mbox{ as Radon measures on }\;\Omega,\;\hbox{ for all }\,k>0$;
\item[$(3)$] $ [{\bf z}, \nu](x)\in \mbox{Sgn}(-u(x))\quad \mathcal{H}^{N-1}-\mbox{ a.e. }
\;x\in \partial\Omega;$
\end{itemize}
\end{definition}

\begin{remark}
Let us explicitly remark that any term in Definition \ref{DefSol} (1) is well defined and that (3) is  condition  is a nowdays standard way to give meaning to the homogeneous Dirichlet boundary datum. It is well known, in fact,  that $BV$ solutions to problems involving the $1$-Laplace operator do not necessarily assume the boundary datum pointwise. Finally, (2) is the weak sense in which  one interprets  ${\bf z}$ as the singular ratio $\frac{D u}{|D u|}$. Here we further weaken this request by truncating as we do not know in general the solutions of to be bounded (see also Section  \ref{r61} below). 
\end{remark}

Our main result is the following
\begin{theorem} \label{exi} Assume \eqref{H}. Then there exists a solution of Problem \eqref{P} in the sense of Definition \ref{DefSol}.
\end{theorem}
In order to prove this result, we consider the approximate solutions of 
\begin{equation}
\label{P1}
\left\{
\begin{array}{rclr}
-\Delta_p u&=&\frac{\lambda}{|x|^p}|u|^{p-2}u+f&\mbox{ in  }  \Omega,\\
u&=&0&\mbox{ on }\partial\Omega,
\end{array}
\right.
\end{equation}
aiming to send  $p\to1^+$ by means of suitable  apriori  estimates independent of  $p$ and then trying to  pass the limit as $p\to 1^+$.
The key  difficulty here relies on the  the influence of zero order perturbation term
$$
\frac{\lambda}{|x|^p}|u|^{p-2}u.
$$
In order to get rid of this drawback one needs to carefully  use Hardy inequality \eqref{HaI} as well as a suitable   truncation argument. 

Furthermore, we will study the asymptotic behaviour of the solutions to \eqref{P1} as $p$ go to $1^+$ considering that \eqref{H} just satisfies the strict inequality (non-extreme case) or the identity (extreme case). More precisely, we shall show the following result:
\begin{theorem}
\label{AsympB}
Assume that $f$ and $\lambda$ are   such that 
$$
\|f\|_{L^N(\Omega)}S_N+\frac{\lambda}{N-1}< 1.
$$
Then
$$
u_p\to 0\qquad \mbox{ a.e. on }\quad\Omega\quad\mbox{ as }\quad p\to 1^+. 
$$
On the other hand, if  $f$ and $\lambda$ are   such that 
$$
\|f\|_{L^N(\Omega)}S_N+\frac{\lambda}{N-1}= 1,
$$
then, there exist $u\in BV(\Omega)$ such that 
$$
u_p\to u\qquad\mbox{ in }\quad L^s(\Omega),\; s\in [1,1^*)\quad\mbox{ as } \quad p\to 1^+.
$$
\end{theorem}
\begin{remark}
Let us stress that, as shown in the proof of Theorem \ref{exi},  the limits of the previous theorem will be solutions to \eqref{P} in the sense of the Definition \ref{DefSol}.
\end{remark}

It is worth mentioning that  although the solutions $u_p$ are, in general,  unbounded (see \cite{OrsinaOrtiz2024, BoccardoOrsinaPeral2006}), their limit can be 0 as established in the previous result. \bk Moreover, the existence of solutions to \eqref{P1} depends only on the restriction of $\lambda$ and the regularity of $f$ (see \cite{OrsinaOrtiz2024, AzoreroAlonso1998, BoccardoOrsinaPeral2006}), however, for \eqref{P}, it will also depend on the 'smallness' of the norm of $f$. This last feature is typical of non-homogeneous problems involving the 1-Laplacian.  For instance, in the torsion problem with  datum $f=1$, this is reduced to the fact that the domain is small (see \cite{K}), while  when the datum belongs to $L^N (\Omega)$, $L^{N,\infty}(\Omega)$ or $W^{-1,\infty}(\Omega)$, this amounts the request of, respectively, small norm (see \cite{ct, mst}). 
\section{ Asymptotic behaviour of the approximating  $p-$Laplacian problems}
Let $1< p<N$ and $0<\beta<((N-p)/p)^p$. Consider the problem
\begin{equation}
\label{Pp}
\left\{
\begin{array}{rclr}
-\Delta_p u_p&=&\frac{\beta}{|x|^p}|u_p|^{p-2}u_p+f&\mbox{ in  }  \Omega,\\
u_p&=&0&\mbox{ on }\partial\Omega,
\end{array}
\right.
\end{equation}
where $f\in W^{-1, p'}(\Omega)$. It is  known that the above problem  has a weak solution $u_p\in W^{1,p}_0(\Omega)$ (see for  instance \cite{AzoreroAlonso1998}), that is, $u_p$ satisfies 
$$
\label{WSolD}
\int_{\Omega}|\nabla u_p|^{p-2}\nabla u_p\cdot\nabla \varphi\,dx=\beta\int_{\Omega}\frac{1}{|x|^p}|u_p|^{p-2}u_p\varphi\,dx+\langle f, \varphi\rangle_{W^{-1,p}(\Omega), W^{1,p}_0(\Omega)}\quad\forall\,\varphi\in W^{1,p}_0(\Omega).
$$

\subsection{Existence for the approximating problems}
For sake completeness and to carefully  keep track of the constants in the a priori estimates, we give a prove of the existence of a solution to \eqref{Pp} by re-adapting the  truncation approach given in \cite{AzoreroAlonso1998}.

\begin{theorem}\label{exip} Assume that $f\in W^{-1, p'}(\Omega)$ .  Then there exists a weak solution $u_p$ of Problem \eqref{Pp}. Moreover, this solution satisfies
$$
\label{EstP}
\|u_p\|_{W^{1,p}_0}\leq \left[\frac{\|f\|_{W^{-1,p}(\Omega)}}{1-\beta\left(\frac{p}{N-p}\right)^p}\right]^\frac{1}{p-1}.
$$
\end{theorem}

\dem
For each $n\in \mathbb{N}$, consider the problem
\begin{equation}
\label{np}
\left\{
\begin{array}{rclr}
-\Delta_p u_n &=&\beta W_n(x)|u_n|^{p-2}u_n+f&\mbox{ in  }  \Omega,\\
u_n&=&0&\mbox{ on }\partial\Omega,
\end{array}
\right.
\end{equation}
where $W_n(x)=\min\{1/|x|^p,n\}$, $n\in \mathbb{N}$. 
In order to find solutions to \eqref{np}, we look at the following minimization problem 
\begin{equation}
\label{MzP}
I_n=\inf_{u\in W^{1,p}_0(\Omega)}\left\{\frac{1}{p}\int_{\Omega}|\nabla u|^p\,dx-\frac{\beta}{p}\int_{\Omega}W_n(x)|u|^p\,dx-\langle f, u\rangle_{W^{-1,p'}(\Omega), W^{1,p}_0(\Omega)}\right\}.
\end{equation}

Let $(u_{n, m})_{n,m\in \mathbb{N}}\subset W^{1,p}_0(\Omega)$ be a minimizing  sequence to \eqref{MzP}, that is,
$$
\frac{1}{p}\int_{\Omega}|\nabla u_{n,m}|^p\,dx-\frac{\beta}{p}\int_{\Omega}W_n(x)|u_{n, m}|^p\,dx-\langle f, u_{n, m}\rangle_{W^{-1,p'}(\Omega), W^{1,p}_0(\Omega)}\to I_n\qquad\mbox{ as }\quad m\to \infty.
$$
Using Hardy inequality, we have that
\begin{align*}
\label{Est2}
&\frac{1}{p}\int_{\Omega}|\nabla u_{n,m}|^p\,dx-\frac{\beta}{p}\int_{\Omega}W_n(x)|u_{n,m}|^p\,dx-\langle f, u_{n, m}\rangle_{W^{-1,p'}(\Omega), W^{1,p}_0(\Omega)}\\
&\geq\frac{1}{p}\left[1-\beta\left(\frac{p}{N-p}\right)^p\right]\|u_{n, m}\|_{W^{1,p}_0(\Omega)}^p-\|f\|_{W^{-1,p}(\Omega)}\|u_{n, m}\|_{W^{1,p}_0(\Omega)}.
\end{align*}
Hence, one can infer that $(u_{n,m})_{n,m\in \mathbb{N}}$ is a bounded sequence in $W^{1,p}_0(\Omega)$ and so, up least to subsequences,
\begin{eqnarray}
\label{CUm}
u_{n, m}\rightharpoonup u_n&\mbox{ in }\quad W^{1,p}_0(\Omega),\nonumber\\
u_{n, m}\rightharpoonup u_n&\mbox{ in }\quad  L^{p^*}(\Omega),\\
u_{n, m}\to u_n&\mbox{ in }\quad L^s(\Omega)\;, s\in[p, p^*),\nonumber
\end{eqnarray}
as $m\to\infty$. 

Using \eqref{CUm} and the weak lower semicontinuity of the energy term, we have that
\begin{align*}
I_n&=\lim_{m\to\infty}\left(\frac{1}{p}\int_{\Omega}|\nabla u_{n,m}|^p\,dx-\frac{\beta}{p}\int_{\Omega}W_n(x)|u_{n, m}|^p\,dx-\langle f, u_{n, m}\rangle_{W^{-1,p'}(\Omega), W^{1,p}_0(\Omega)}\right)\\
&\geq \frac{1}{p}\int_{\Omega}|\nabla u_n|^p\,dx-\frac{\beta}{p}\int_{\Omega}W_n(x)|u_n|^p\,dx-\langle f, u_n\rangle_{W^{-1,p'}(\Omega), W^{1,p}_0(\Omega)}\\
&\geq I_n,
\end{align*}
which implies that, the infimum in \eqref{MzP} is attained. Consequently, since 
$$
u\mapsto \frac{1}{p}\int_{\Omega}|\nabla u|^p\,dx-\frac{\beta}{p}\int_{\Omega}W_n(x)|u|^p\,dx-\langle f, u\rangle_{W^{-1,p}(\Omega), W^{1,p}(\Omega)} 
$$
is differentiable in $ W^{1,p}_0(\Omega)$,  the weak form of the Euler equation holds, that is 
\begin{equation}
\label{WUn}
\int_{\Omega}|\nabla u_n|^{p-2}\nabla u_n\cdot\nabla \varphi\,dx-\beta\int_{\Omega}W_n(x)|u_n|^{p-2}u_n\varphi\,dx-\langle f, \varphi\rangle_{W^{-1, p'}(\Omega), W^{1,p}_0(\Omega)}=0\quad\forall\;\varphi\in W^{1,p}_0(\Omega).
\end{equation}

Use $u_n$ as test function in \eqref{WUn}. Then, by Hardy and Sobolev inequalities,
\begin{equation}
\label{BUnifN}
\|u_n\|_{W^{1,p}_0(\Omega)}\leq \left[\frac{\|f\|_{W^{-1,p}(\Omega)}}{1-\beta\left(\frac{p}{N-p}\right)^p}\right]^\frac{1}{p-1}\quad\mbox{ for all }\; n\in \mathbb{N}.
\end{equation}
Consequently, up to subsequences,
\begin{align*}
\label{WCUn}
u_n\rightharpoonup u_p&\mbox{ in }\;W^{1,p}_0(\Omega)\\
 u_n\to u_p& \mbox{ in }\;L^s(\Omega),\; s\in [p, p^*),
\end{align*}
as $n\to \infty$. 

{In order to pass to the limit in \eqref{WUn} as $n\to \infty$, first we will show that
\begin{equation}
\label{ConvL1}
W_n(x)|u_n|^{p-2}u_n\to \frac{1}{|x|^p}|u_p|^{p-2}u_p\quad\mbox{ in }\; L^1(\Omega)\quad\hbox{ as }\; n\to \infty.
\end{equation}
Since $u_n\to u_p$ in $L^s(\Omega)$, $s\in [p, p^*)$, by Vainberg theorem, up least to subsequences,
$$
u_n(x)\to u_p(x)\quad\mbox{ a.e. }\; x\in\Omega,\quad\hbox{ as }\;n\to \infty
$$
and there exists $g\in L^s(\Omega)$ such that
$$
|u_n(x)|\leq g(x)\quad\mbox{ a.e. }\; x\in \Omega.
$$
Thus
\begin{equation}
\label{PCWn}
W_n(x)|u_n(x)|^{p-2}u_n(x)\to \frac{1}{|x|^p}|u_p(x)|^{p-2}u_p(x)\quad\hbox{ a.e. } \;x\in \Omega\quad\hbox{ as}\; n\to\infty
\end{equation}
and
$$
|W_n(x)|u_n(x)|^{p-2}u_n(x)|\leq \frac{1}{|x|^p} g(x)^{p-1}\quad\hbox{ a.e. }\; x\in \Omega.
$$

{\bf Claim 1:}
$$
\int_{\Omega} \frac{1}{|x|^p} g(x)^{p-1}\,dx<\infty.
$$
Indeed, consider $1<r< p^*/(p-1)$ and so $p^*/(p^*-p+1)<r'=\frac{r}{r-1}$. We may choose
$$
r'=\frac{p^*(1+\epsilon)}{p^*-p+1},
$$
where $\epsilon>0$ is such that $r'<N/p$. Thus, by H\"older inequality, 
$$
\int_{\Omega}\frac{1}{|x|^p}g(x)^{p-1}\,dx\leq \left(\int_{\Omega}\frac{1}{|x|^{pr'}}\,dx\right)^\frac{1}{r'}\left( \int_{\Omega} g(x)^{(p-1)r}\,dx \right)^\frac{1}{r}<\infty.
$$

Consequently, from \eqref{PCWn} and Claim 1, by Lebesgue's dominated convergence theorem, \eqref{ConvL1} holds.

{\bf Claim 2: }$\nabla u_n\to \nabla u_p$ in $L^q(\Omega)$ for all $q<p$.

We can reason as in the proof of  \cite[Theorem 2.1]{BoccardoMurat1992} to get 
\begin{equation}
\label{Gun}
\nabla u_n(x)\to \nabla u_p(x)\quad{ a.e. }\quad x\in\Omega\quad\hbox{ as }\quad n\to \infty.
\end{equation}

Now, for each $q<p$ notice that
\begin{equation}
\label{Eq}
\int_{E}|\nabla u_n|^q\,dx<\epsilon\quad\forall\, n\in \mathbb{N}\quad\hbox{ whenever} \quad|E|<\delta,
\end{equation}
where
$$
\delta<\left\{\epsilon\left[\frac{1-\beta\left(\frac{p}{N-p}\right)^p}{\|f\|_{W^{-1,p}(\Omega)}}\right]^\frac{q}{p-1}\right\}^\frac{p}{p-q}.
$$
Then, from \eqref{Gun} and \eqref{Eq}, by Vitali's theorem, Claim 2 holds. Thus, letting $n\to \infty$ in \eqref{WUn}, from \eqref{ConvL1}, we obtain that
$$
\label{Wu}
\int_{\Omega}|\nabla u_p|^{p-2}\nabla u_p\cdot\nabla \varphi\,dx=\beta\int_{\Omega}\frac{1}{|x|^p}|u_p|^{p-2}u_p\varphi\,dx+\langle f, \varphi\rangle_{W^{-1,p}(\Omega), W^{1,p}_0(\Omega)}\quad \mbox{for all }\;\varphi\in C^\infty_c(\Omega).
$$
Moreover, by a density argument, it is true for all $\varphi\in W^{1,p}_0(\Omega)$. 

On the other hand, passing to the limit in \eqref{BUnifN}, we have that $u_p$ satisfies
$$
\label{BUnifP}
\|u_p\|_{W^{1,p}_0(\Omega)}\leq \left[\frac{\|f\|_{W^{-1,p}(\Omega)}}{1-\beta\left(\frac{p}{N-p}\right)^p}\right]^\frac{1}{p-1}.
$$
\fim

\subsection{Basic estimates}

Recalling that $\lambda< N-1$ then there exists $\bar p\in (1, N)$ such that
$$
\label{EstLbd}
\lambda<\left(\frac{N-p}{p}\right)^p\qquad \hbox{ for all }\; p\in (1,\bar p].
$$

For each $p\in (1, \bar p]$ and $n\in \mathbb{N}$ consider the problem
$$
\label{Lnp}
\left\{
\begin{array}{rclr}
-\Delta_p u &=&\frac{\lambda}{|x|^p} |u|^{p-2}u+f&\mbox{ in  }  \Omega,\\
u&=&0&\mbox{ on }\partial\Omega,
\end{array}
\right.
$$
where $f\in L^N(\Omega)$. Since $f\in L^N(\Omega)\subset W^{-1, p'}(\Omega)$, by Theorem \ref{exip}, there exists a weak solution $u_p$ to this problem, that is,
\begin{equation}
\label{WSol}
\int_{\Omega}|\nabla u_p|^{p-2}\nabla u_p\cdot\nabla \varphi\,dx=\lambda\int_{\Omega}\frac{1}{|x|^p}|u_p|^{p-2}u_p\varphi\,dx+\int_{\Omega}f\varphi\,dx\quad\forall\,\varphi\in W^{1,p}_0(\Omega).
\end{equation}
Moreover, the family of solutions $(u_p)_{p>1}$ satisfies the following inequality
\begin{equation}
\label{B}
\|u_p\|_{W^{1,p}_0(\Omega)}\leq \left[\frac{\|f\|_{L^N(\Omega)}S_N}{1-\lambda\left(\frac{p}{N-p}\right)^p}\right]^\frac{1}{p-1}|\Omega|^\frac{1}{p}\qquad\mbox{ for all }\; p\in(1,\bar p].
\end{equation}

 It follows from the inequality above and condition \eqref{H} the following result.
 
\begin{lemma}
\label{UB}
 Assume that $f\in L^{N}(\Omega)$ is such that   \eqref{H} holds.  Then the family of solutions $(u_p)_{p>1}$  is uniformly bounded in $BV(\Omega)$. Furthermore, there exists $u \in BV (\Omega)$ such that, up least to subsequences,
 \begin{eqnarray*}
 \label{ConLp}
 u_p\to u&\mbox{in }\quad L^s(\Omega),\; s\in [1, 1^*),\\
 \nonumber u_p \rightharpoonup u&\mbox{ in } \quad L^{1^*}(\Omega),\\
 \nonumber \nabla u_p\rightharpoonup Du&*-\hbox{weak}\quad\hbox{as measures in}\quad\Omega,
 \end{eqnarray*}
 as $ p\to 1^+$.
\end{lemma}

\dem
 From \eqref{B}, by applying \eqref{H}, we have that
 \begin{equation}
 \label{EstUp0}
 \int_{\Omega}|\nabla u_p|^p\,dx\leq \left[\frac{1-\frac{\lambda}{N-1}}{1-\lambda\left(\frac{p}{N-p}\right)^p} \right]^\frac{p}{p-1}|\Omega| \quad \mbox{ for all }\quad p\in(1,\bar p].
 \end{equation}
Since
 $$
 \lim_{p\to1^+}\left[\frac{1-\frac{\lambda}{N-1}}{1-\lambda\left(\frac{p}{N-p}\right)^p} \right]^\frac{p}{p-1}=e^{-\frac{\lambda[N-(N-1)\ln(N-1)]}{(N-1)(N-1-\lambda)}},
 $$
 there exists $C>0$ such that
 \begin{equation}
 \label{BC}
 \int_{\Omega}|\nabla u_p|^p\,dx\leq C|\Omega|\qquad\mbox{ for all }\; p\in (1,\bar p].
 \end{equation}
 By Young's inequality and the fact that $u_p=0$ on $\partial\Omega$, it follows that
$$
 \label{BU}
 \int_{\Omega}|\nabla u_p|\,dx+\int_{\partial\Omega}|u_p|\,dx\leq\frac{1}{p}\int_{\Omega}|\nabla u_p|^p\,dx+\frac{p-1}{p}|\Omega|\leq \frac{1}{p}C|\Omega|+\frac{p-1}{p}|\Omega|\leq (C+1)|\Omega|,
$$
for all $p\in (1,\bar p]$. But this, in turn, implies that $(u_p)_{p>1}$ is a uniformly bounded in $BV(\Omega)$. 
Consequently, by the compact embedding of $BV(\Omega)$, there exists $u\in BV(\Omega)$ such that, up least to subsequences,
\begin{equation}
\label{ConvUp}
u_p\to u\quad\mbox{in}\quad L^s(\Omega),\; s\in [1, 1^*) \quad\mbox{ as }\; p\to 1^+,
\end{equation}
 by Sobolev's inequality, 
$$
u_p\rightharpoonup u\quad\mbox{ in }\;L^{1^*}(\Omega)\quad\mbox{ as }\; p\to 1^+
$$
and 
$$
\nabla u_p\rightharpoonup Du\qquad*-\hbox{weak}\quad\hbox{as measures in }\quad\Omega.
$$
\fim

\begin{lemma}
\label{ConvZ}
Let $q\geq 1$. There exists ${\bf z}\in L^\infty(\Omega,\mathbb{R}^N)$ with $\|{\bf z}\|_{L^\infty(\Omega)}\leq 1$ such that 
$$
|\nabla u_p|^{p-2}\nabla u_p\rightharpoonup{\bf z}\quad \mbox{ in }\;L^q(\Omega,\mathbb{R}^N)\quad\mbox{ as }\;p\to 1^+.
$$
\end{lemma}			

\dem Consider $q=1$. We will show that $(|\nabla u_p|^{p-2}\nabla u_p)_{p>1}$ is equi-integrable in $L^1(\Omega,\mathbb{R}^N)$. Indeed, by \eqref{BC} and the H\"older's inequality, it follows that 
$$
\int_{\Omega}|\nabla u_p|^{p-1}\,dx\leq \left(\int_{\Omega}|\nabla u_p|^p\,dx\right)^\frac{p-1}{p}|\Omega|^\frac{1}{p}\leq (C+1)|\Omega|\quad\hbox{ for all }\;p\in (1, \bar p].
$$
Moreover, for any measurable subset $E\subset \Omega$, we have 
$$
\left|\int_{E}|\nabla u_p|^{p-2}\nabla u_p\,dx\right|\leq \int_{E}|\nabla u_p|^{p-1}\,dx\leq \left(\int_{E}|\nabla u_p|^p\,dx\right)^\frac{p-1}{p}|E|^\frac{1}{p}\quad\hbox{ for all }\; p\in (1, \bar p].
$$
Hence $(|\nabla u_p|^{p-2}\nabla u_p)_{p>1}$ is equi-integrable in $L^1(\Omega,\mathbb{R}^N)$ and so,  by Pettis Theorem,  it is weakly relatively compact in $L^1(\Omega,\mathbb{R}^N)$. Consequently, there exists ${\bf z}\in L^1(\Omega,\mathbb{R}^N)$ such that, up  to subsequences, 
$$
|\nabla u_p|^{p-2}\nabla u_p\rightharpoonup {\bf z}\quad\mbox{ in } L^1(\Omega,\mathbb{R}^N) \quad \mbox{ as } \; p\to 1^+.
$$

On the other hand, for each $q>1$, we have that
\begin{equation}\label{idenz}
\int_{\Omega}||\nabla u_p|^{p-2}\nabla u_p|^q\,dx \leq \left(\int_{\Omega}|\nabla u_p|^p\,dx\right)^\frac{(p-1)q}{p}|\Omega|^\frac{p-(p-1)q}{p}\leq (C+1)|\Omega|\quad\hbox{ for all } \; p\in (1, \bar p].
\end{equation}
Then there exists ${\bf z}_q\in L^q(\Omega,\mathbb{R}^N)$ such that, up to  subsequences,
$$
|\nabla u_p|^{p-2}\nabla u_p\rightharpoonup {\bf z}_{q}\quad \mbox{ in } L^q(\Omega,\mathbb{R}^N)\quad\mbox{ as }\; p\to 1^ +.
$$
By a diagonal argument,  we can find ${\bf z}$ independent of $q$ such that it also be the limit of $(|\nabla u_p|^{p-2}\nabla u_p)_{p>1}$ in $L^q(\Omega, \mathbb{R}^N)$. Thus, letting $p\to 1^+$ in \eqref{idenz}, by lower semicontinuity of norm in $L^q(\Omega, \mathbb{R}^N)$, we infer that
$$
\|{\bf z}\|_{L^q(\Omega,\mathbb{R}^N)}\leq \liminf_{p\to 1^+}\||\nabla u_p|^{p-2}\nabla u_p\|_{L^q(\Omega)}\leq |\Omega|^\frac{1}{q}.
$$
But this, in turn, letting $q\to\infty$, implies that 
$$
\|{\bf z}\|_{L^\infty(\Omega,\mathbb{R}^N)}\leq 1.
$$
\fim

\section{Proof of Theorem 3.1}
In this section, we shall show that limit in \eqref{ConvUp}, $u$, with ${\bf z}$ as in \eqref{ConvZ}, is a solution for \eqref{P} in sense of the Definition \ref{DefSol}. 

First, we denote 
$$
\{u\neq 0 \}=\{ x \in \Omega\>:\> u(x)\neq 0 \;\hbox{ a.e.}\}\quad\hbox{ and }\quad\{u = 0\} = \{ x\in \Omega\>;\> u(x) = 0\;\hbox{ a.e.}\}.
$$

{\bf  (1)} We shall prove that 
\begin{equation}   
\label{ConvS1}
\int_{\Omega}{\bf z}\cdot\nabla \varphi\,dx=\lambda\int_{\Omega}\frac{1}{|x|}s(x)\varphi\,dx+\int_{\Omega} f\varphi\,dx \quad \mbox{ for all } \; \varphi\in C^\infty_c(\Omega),
\end{equation}
with a   measurable selection $s(x)\in \hbox{Sgn}(u(x))$,  $\hbox{Sgn}(\cdot)$  defined in \eqref{Sgn}. 

In order to prove \eqref{ConvS1}, we will pass to the limit in \eqref{WSol}. Notice that, by the Lemma \ref{ConvZ}, we get that
\begin{equation}
\label{1ConvZ}
\lim_{p\to1^+}\int_{\Omega}|\nabla u_p|^{p-2}\nabla u_p\cdot\nabla \varphi\,dx=\int_{\Omega}{\bf z}\cdot\nabla \varphi\,dx.
\end{equation}
{\bf Claim:}
$$
\label{Conv0}
\lim_{p\to1^+}\int_{\Omega}\frac{1}{|x|^p}|u_p|^{p-2}u_p\varphi\,dx=\int_{\Omega}\frac{1}{|x|}s(x)\varphi\,dx.
$$
where $s(x)\in \hbox{Sgn}(u(x))$ a.e. $x\in \Omega$. We write
$$
\int_{\Omega}\frac{1}{|x|^p}|u_p|^{p-2}u_p\varphi\,dx=\int_{ \{u\neq 0\} }\frac{1}{|x|^p}|u_p|^{p-2}u_p\varphi\,dx + \int_{\{u = 0\}}\frac{1}{|x|^p}|u_p|^{p-2}u_p\varphi\,dx.
$$
Since $u_p\to u$ in $L^s(\Omega)$, $s\in [1,1^*)$, by Vainberg theorem, up to subsequences, 
\begin{equation}
\label{ConvPu}
u_p(x)\to u(x)\quad\mbox{ a.e. }\quad x\in \Omega \quad \mbox{ as } \; p\to 1^+
\end{equation}
and there exists $g\in L^s(\Omega)$, $s\in [1, 1^*)$, such that
\begin{equation}
\label{g}
|u_p(x)|\leq g(x)\quad\mbox{ a.e. } \; x\in \Omega.
\end{equation}
Thus
$$
\frac{1}{|x|^p}|u_p(x)|^{p-2}u_p(x)\varphi(x)\to \frac{1}{|x|}\frac{u(x)}{|u(x)|}\varphi(x)\quad\mbox{ a.e. }\;x\in \{u \neq 0\}\quad\mbox{as} \quad p\to1^+
$$
and, using the Young inequality, 
$$
\left|\frac{1}{|x|^p}|u_p(x)|^{p-2}u_p(x)\varphi(x)\right|\leq \|\varphi\|_{L^\infty(\Omega)}\max\left\{\frac{1}{|x|^{\bar p}},1\right\}\left(|g(x)|^{\bar p-1}+1\right)\quad\hbox{ a.e. }\;x\in \;\forall \, p\in (1, \bar p].
$$
It is not difficult to show that the term on right-side hand above belongs to $L^1(\Omega)$. Indeed, using H\"older inequality with exponents $r\in(1, 1^*/(\bar p-1))$ and $r'=1^*(1+\epsilon)/(1^*-p+1)$ for some $\epsilon>0$ such that $r'<N/\bar p$, we obtain that
\begin{equation}
\label{L1}
\int_{\Omega}\max\left\{\frac{1}{|x|^{\bar p}},1\right\}\left(|g(x)|^{\bar p-1}+1\right)\,dx\leq \left(\int_{\Omega}\max\left\{\frac{1}{|x|^{\bar p}},1\right\}^{r'}\,dx\right)^\frac{1}{r'}\left(\int_{\Omega}(g(x)^{\bar p-1}+1)^r\,dx\right)^\frac{1}{r}<\infty.
\end{equation}
Consequently, by Lebesgue's dominated convergence theorem, 
\begin{equation}
\label{UN0}
\lim_{p\to 1^+}\int_{\{u \neq 0\}}\frac{1}{|x|^p}|u_p|^{p-2}u_p\varphi\,dx=\int_{\{u\neq 0\}}\frac{1}{|x|}\frac{u}{|u|}\varphi\,dx.
\end{equation}

On the other hand, for each $q>1$, we have that
\begin{align}
\label{Estq}
\nonumber\int_{\{u=0\}}||u_p|^{p-2}u_p|^q\,dx&\leq \left(\int_{\{u = 0\}}|u_p|^p\,dx\right)^\frac{(p-1)q}{p}|\Omega|^{\frac{p-(p-1)q}{p}}\\
&\leq S_{N,p}^{(p-1)q}\left(\int_{\Omega}|\nabla u_p|^p\,dx\right)^\frac{(p-1)q}{p}|\Omega|^\frac{p-(p-1)q}{p}\\
\nonumber&\leq S_{N,p}^{(p-1)q}C^\frac{(p-1)q}{p}|\Omega|,
\end{align}
for all $p\in (1, \bar{p}]$. Hence $(|u_p|^{p-2}u_p)_{p>1}$ is a bounded sequence in $L^q(\{u = 0\})$ and so, there exists $v_q\in L^q(\{u=0\})$ such that, up least to subsequences,
$$
|u_p|^{p-2}u_p\rightharpoonup v_q\quad\mbox{ in }\; L^q(\{u=0\})\quad\hbox{ as }\; p\to 1^+.
$$
Again by a diagonal argument, there exists $v$, independent of $q>1$, such that
\begin{equation}
\label{WCq}
|u_p|^{p-2}u_p\rightharpoonup v\quad\mbox{ in }\; L^q(\{u=0\})\;\mbox{ for all }\;q>1\quad\hbox{ as }\quad p\to 1^+.
\end{equation}
But this, in turn, setting $p\to 1^+$ in \eqref{Estq}, implies that
$$
\|v\|_{L^q(\{u =0\})}\leq |\Omega|^\frac{1}{q}\quad\mbox{ for all }\; q>1.
$$
Passing to the limit as $q\to \infty$ above, we obtain
$$
\|v\|_{L^\infty(\{u = 0\})}\leq 1.
$$
Consequently, using \eqref{WCq}, 
\begin{equation}
\label{U0}
\lim_{p\to 1^+}\int_{\{u=0\}}\frac{1}{|x|^p}|u_p|^{p-2}u_p\varphi\,dx=\int_{\{u = 0\}}\frac{1}{|x|}v\varphi\,dx.
\end{equation}

Thus, from \eqref{UN0} and \eqref{U0}, we get that
$$
\label{CSgn}
\lim_{p\to 1^+}\int_{\Omega}\frac{1}{|x|^p}|u_p|^{p-2}u_p\varphi\,dx=\int_{\Omega}\frac{1}{|x|}s(x)\varphi\,dx,
$$
where
$$
s(x)=\left\{
\begin{array}{lcr}
\frac{u(x)}{|u(x)|}& \hbox{if} &u(x)\neq 0\\
v(x) &\hbox{if }& u(x)=0
\end{array}
\right\}\in Sgn(u(x))\quad\hbox{ a.e. }\;\;x\in \Omega,
$$
and so, the Claim is verified.

Therefore, setting $p\to 1^+$ in \eqref{WSol}, by \eqref{1ConvZ}, (1) of Definition \ref{DefSol} holds.

{\bf (2)} It follows, by Anzellotti's theory, that
\begin{equation}\label{rev}
({\bf z}, DT_k(u))\leq |DT_k(u)|\quad\mbox{ as Radon measures on }\;\Omega.
\end{equation} 
In order to prove the reverse inequality we use $T_k(u_p)\varphi$ as test function in \eqref{WSol}. 
Then
 $$
 \begin{aligned}
 \int_{\Omega}\varphi|\nabla T_k(u_p)|^p\,dx+\int_{\Omega}T_k(u_p)|\nabla u_p|^{p-2}\nabla u_p\cdot\nabla \varphi\,dx \\ =\lambda\int_{\Omega}\frac{1}{|x|^p}|u_p|^{p-2}u_pT_k(u_p)\varphi\,dx+\int_{\Omega}fT_k(u_p)\varphi\,dx.
 \end{aligned}
 $$
By Young's inequality, we have that
 \begin{align}
 \nonumber\int_{\Omega}\varphi|DT_k(u_p)|\,dx&\leq\frac{1}{p}\int_{\Omega}\varphi|\nabla T_k(u_p)|^p\,dx+\frac{p-1}{p}\int_{\Omega}\varphi\,dx\\
 \label{LowSem}
 &=\lambda\int_{\Omega}\frac{1}{|x|^p}|u_p|^{p-2}u_pT_k(u_p)\varphi\,dx+\int_{\Omega}fT_k(u_p)\varphi\,dx\\
 \nonumber&\qquad-\int_{\Omega}T_k(u_p)|\nabla u_p|^{p-2}\nabla u_p\cdot \nabla\varphi\,dx+\frac{p-1}{p}\int_{\Omega}\varphi\,dx.
 \end{align}
It is not difficult to show that
\begin{equation}
\label{ConvS2}
\lim_{p\to1^+}\int_{\Omega}\frac{1}{|x|^p}|u_p|^{p-2}u_pT_k(u_p)\varphi\,dx=\int_{\Omega}\frac{1}{|x|} s(x)T_k(u)\varphi\,dx,
\end{equation}
with a measurable selection $s(x)\in \hbox{Sgn}(u(x))$. Indeed, by \eqref{ConvPu} and  \eqref{g}, we have
$$
\frac{1}{|x|^p}|u_p(x)|^{p-2}u_p(x) T_k(u_p(x))\varphi(x)\to \frac{1}{|x|}\frac{u(x)}{|u(x)|}T_k(u(x))\varphi(x)\quad\mbox{ a.e. }\; x\in \{ u\neq 0 \}\quad\mbox{ as }\;p\to 1^+
$$
and
$$
\left|\frac{1}{|x|^p}|u_p(x)|^{p-2}u_p(x) T_k(u_p(x))\varphi(x)\right|\leq k\|\varphi\|_{L^\infty(\Omega)}\max\left\{\frac{1}{|x|^{\bar p}},1\right\}\left(\frac{p-1}{\bar p-1}g(x)^{\bar p-1}+\frac{\bar p-p}{\bar p-1}\right),
$$
 for all $p\in (1, \bar p]$. By same arguments in \eqref{L1}, the term on right side-hand above belongs to $L^1(\Omega)$. Thus
$$
\lim_{p\to 1^+}\int_{ \{u \neq 0\}}\frac{1}{|x|^p}|u_p|^{p-2}u_pT_k(u_p)\varphi\,dx=\int_{\{ u\neq 0\}}\frac{1}{|x|}\frac{u}{|u|}T_k(u)\varphi\,dx.
$$
Moreover, by same argument in \eqref{U0}, there exists $v\in L^\infty(\Omega)$ such that $\|v\|_{L^\infty(\Omega)}\leq 1$ and 
$$
\lim_{p\to 1^+}\int_{\{u = 0\}}\frac{1}{|x|^p}|u_p|^{p-2}u_pT_k(u_p)\varphi\,dx=\int_{\{u= 0\}}\frac{1}{|x|}vT_k(u)\varphi\,dx.
$$
Consequently, \eqref{ConvS2} holds.

One can also show that
\begin{equation}
\label{CTk}
\lim_{p\to1^+}\int_{\Omega}T_k(u_p)|\nabla u_p|^{p-2}\nabla u_p\cdot\nabla \varphi\,dx=\int_{\Omega}T_k(u){\bf z}\cdot\nabla \varphi\,dx.
\end{equation}
We write
\begin{align*}
&\int_{\Omega}T_k(u_p)|\nabla u_p|^{p-2}\nabla u_p\cdot\nabla \varphi\,dx-\int_{\Omega}T_k(u){\bf z}\cdot\nabla\varphi\,dx\\
&=\int_{\Omega}(T_k(u_p)-T_k(u))|\nabla u_p|^{p-2}\nabla u_p\cdot\nabla \varphi\,dx+\int_{\Omega}T_k(u)(|\nabla u_p|^{p-2}u_p-{\bf z})\cdot \nabla\varphi\,dx.
\end{align*}
Fix $q\in(\bar p,1^*)$. Using H\"older's inequality and \eqref{BC} we obtain that
\begin{equation*}\begin{aligned}
&&\left|\int_{\Omega}(T_k(u_p)-T_k(u))|\nabla u_p|^{p-2}\nabla u_p\cdot\nabla\varphi\,dx\right|  &\leq& \||\nabla \varphi|\|_{L^\infty(\Omega)}\|u_p\|_{W^{1,p}_0(\Omega)}^{p-1}\|T_k(u_p)-T_k(u)\|_{L^q(\Omega)}|\Omega|^{\frac{1}{p}-\frac{1}{q}}.\end{aligned}
\end{equation*}
Since, as an easy consequence of \eqref{ConvUp}, $T_k(u_p)\to T_k(u)$ in $L^q(\Omega)$, we have that
\begin{equation}
\label{GConvUp}
\lim_{p\to 1^+}\int_{\Omega}(T_k(u_p)-T_k(u))|\nabla u_p|^{p-2}\nabla u_p\cdot\nabla\varphi\,dx=0.
\end{equation}
On the other hand, since $|\nabla u_p|^{p-2}\nabla u_p\rightharpoonup {\bf z}$ in $L^1(\Omega,\mathbb{R}^N)$,
\begin{equation}
\label{GConvWUp}
\lim_{p\to 1^+}\int_{\Omega}T_k(u)(|\nabla u_p|^{p-2}u_p-{\bf z})\cdot \nabla\varphi\,dx=0.
\end{equation}
Thus, from \eqref{GConvUp} and \eqref{GConvWUp}, \eqref{CTk} holds.

Therefore, letting $p\to 1^+$ in \eqref{LowSem}, by the lower semicontinuous of $u\mapsto\int_{\Omega}\varphi|Du|$, \eqref{ConvS2} and \eqref{CTk}, we get that
\begin{align}
\label{EstTku}
\nonumber\int_{\Omega}\varphi|DT_k(u)|& \leq \lambda\int_{\Omega}\frac{1}{|x|}s(x)T_k(u)\varphi\,dx+\int_{\Omega}fT_k(u)\varphi\,dx-\int_{\Omega}T_k (u){\bf z}\cdot\nabla \varphi\,dx\\
&=-\int_{\Omega}(T_k(u))^*\hbox{ div}\,{\bf z}-\int_{\Omega}T_k(u){\bf z}\cdot\nabla \varphi\,dx\\
\nonumber&=\langle ({\bf z}, DT_k(u)), \varphi \rangle.
\end{align}
where $s(x)\in \hbox{Sgn}(u(x))$ a.e. $x\in \Omega$. Hence the reverse inequality of \eqref{rev} holds and so, (2) of Definition \ref{DefSol}.

{\bf (3)} By Proposition \ref{Prop2}, we have that 
$$
\label{T}
|[{\bf z}, \nu](x)|\leq \|{\bf z}\|_{L^\infty(\Omega,\mathbb{R}^N)}\leq 1\quad\mathcal{H}^{N-1}-\hbox{ a.e.}\; x\in \partial\Omega.
$$
Consequently
$$
\label{Inq}
\int_{\partial\Omega}(|u|+ u [{\bf z}, \nu ])\,d\mathcal{H}^{N-1}\geq 0.
$$

We shall prove that the reverse inequality holds. Use $T_k(u_p)$ as test function in \eqref{WSol}. Then
$$
\int_{\Omega}|\nabla T_k(u_p)|^p\,dx=\lambda\int_{\Omega}\frac{1}{|x|^p}|u_p|^{p-2}u_pT_k(u_p)\,dx+\int_{\Omega}fT_k(u_p)\,dx.
$$
Using Young's inequality, we get that
\begin{align}
\label{EstDUp}
\int_{\Omega}|\nabla T_k(u_p)|\,dx&\leq \frac{1}{p}\int_{\Omega}|\nabla T_k(u_p)|^p\,dx+\frac{p-1}{p}|\Omega|\\
&\leq \lambda\int_{\Omega}\frac{1}{|x|^p}|u_p|^{p-2}u_pT_k(u_p)\,dx+\int_{\Omega}fT_k(u_p)\,dx+\frac{p-1}{p}|\Omega|.\nonumber
\end{align}
Note that, by the same argument to show \eqref{ConvS2}, we have that
$$
\lim_{p\to 1^+}\int_{\Omega}\frac{1}{|x|^p}|u_p|^{p-2}u_pT_k(u_p)\,dx=\int_{\Omega}\frac{1}{|x|}s(x)T_k(u)\,dx, 
$$
where $s(x)\in \hbox{Sgn}\,(u(x))$ a.e. $x\in \Omega$. Moreover, since $T_k(u_p)\to T_k(u)$ in $L^q(\Omega)$, $q>1$,
$$
\lim_{p\to1^+}\int_{\Omega}fT_k(u_p)\,dx=\int_{\Omega}fT_k(u)\,dx.
$$

Setting $p\to 1^+$ in \eqref{EstDUp}, by the lower semicontinuous of $u\mapsto|Du|(\Omega)$ and the fact that $T_k(u_p)\to T_k(u)$ in $L^1(\Omega)$ as $p\to1^+$, we get
\begin{align*}
\int_{\Omega}|DT_k(u)|+\int_{\partial \Omega}|T_k(u)|\,d\mathcal{H}^{N-1}&\leq \lambda\int_{\Omega}\frac{1}{|x|}s(x)T_k(u)\,dx+\int_{\Omega}fT_k(u)\,dx\\
&=-\int_{\Omega}(T_k(u))^*\hbox{div}\,{\bf z}\\
&=\int_{\Omega}({\bf z}, DT_k(u))-\int_{\partial\Omega}[{\bf z}, \nu] T_k(u)\,d\mathcal{H}^{N-1},
\end{align*}
where $s(x)\in \hbox{Sgn}(u(x))$ a.e. $x\in \Omega$. Hence, by (2) of Definition \ref{DefSol}, 
\begin{equation}
\label{CInq}
\int_{\partial\Omega}|T_k(u)|\,d\mathcal{H}^{N-1}+\int_{\partial\Omega} [{\bf z}, \nu] T_k(u)\,d\mathcal{H}^{N-1}\leq 0\quad\hbox{ for all }\; k>0.
\end{equation}
Note that $T_k(u)\to u$ in $L^1(\Omega)$ and $\int_{\Omega}|DT_k(u)|\leq \int_{\Omega}|Du|$ for all $k>0$, and so,
$$
\int_{\Omega}|Du|=\lim_{k\to\infty}\int_{\Omega}|DT_k(u)|,
$$
in other words $T_k(u)$ strictly converges to $u$. Thus
$$
T_k(u)\to u\quad\hbox{ in } L^1(\partial\Omega)\quad\hbox{ as }\; k\to \infty.
$$
Hence, setting $k\to \infty$ in \eqref{CInq}, we obtain that
$$
\int_{\partial\Omega}|u|\,d\mathcal{H}^{N-1}+\int_{\partial\Omega} [{\bf z}, \nu] u\,d\mathcal{H}^{N-1}\leq 0.
$$
Thus the contrary inequality holds and $(3)$ of the Definition \ref{DefSol} is verified.

Therefore $u$ satisfies $(1)-(3)$ of the Definition \ref{DefSol}.

\begin{remark}
\label{uR}
Let us stress that, in general, one is not able to show that $ u^*\in L^1(\Omega, |\hbox{div}\,{\bf z}|)$ and $u^+, u^-\in L^1(J_u, |\hbox{div}\,{\bf z}|)$. If this is the case instead,  passing to the limit in \eqref{EstTku}, and using by Proposition \ref{Tku}, we obtain that
$$
\int_{\Omega}\varphi|Du|\leq \langle ({\bf z}, Du),\varphi\rangle \quad\hbox{ for all nonnegative}\; \varphi\in C^\infty_c(\Omega).
$$
Hence, as $({\bf z}, Du)\leq |Du|$ (again by Proposition \ref{Tku}), it follows that
$$
({\bf z}, Du)=|Du|\quad\hbox{ as Radon measure on }\;\Omega.
$$
\end{remark}

\section{Proof of the Theorem \ref{AsympB} and some first explicit examples}
 
In this section, we shall prove the Theorem \ref{AsympB} by splitting the proof in the two possible (say non-extreme and extreme) cases:  

{\bf Case 1.} Suppose that $\|f\|_{L^N(\Omega)}S_N+ \frac{\lambda}{N-1}<1$. 

There exists $ \hat p>1$ such that
$$
\|f\|_{L^N(\Omega)}S_N+\lambda\left(\frac{p}{N-p}\right)^p<1\qquad\mbox{ for all }\quad p\in (1, \hat p].
$$
Thus
\begin{equation*}
\lim_{p\to 1^+}\left[\frac{\|f\|_{L^N(\Omega)}S_N}{1-\lambda\left(\frac{p}{N-p}\right)^p} \right]^\frac{p}{p-1}=0.
\end{equation*}
But this, in turn, by \eqref{B}, implies that
$$
\lim_{p\to 1^+}\int_{\Omega}|\nabla u_p|^p\,dx=0.
$$
Consequently, by Sobolev and Young inequalities, we have that
$$
\lim_{p\to 1^+}\int_{\Omega}|u_p|^\frac{N}{N-1}\,dx\leq S_N\lim_{p\to 1^+}\int_{\Omega}|\nabla u_p|\leq S_N\lim_{p\to 1^+}\left(\frac{1}{p}\int_{\Omega}|\nabla u_p|^p\,dx+\frac{p-1}{p}|\Omega|\right)=0
$$
and, in particular 
$$
\lim_{p\to1^+}u_p(x)=0\qquad\mbox{ a.e. }\;x\in\Omega.
$$

{\bf Case 2.} Assume that $\|f\|_{L^N(\Omega)}S_N+\frac{\lambda}{N-1}=1$. 

We start by showing that
$$
\lim_{p\to 1^+}\left[\frac{1-\frac{\lambda}{N-1}}{1-\lambda\left(\frac{p}{N-p}\right)^p} \right]^\frac{p}{p-1}<\infty.
$$
Indeed, since 
$$
\left[\frac{1-\frac{\lambda}{N-1}}{1-\lambda\left(\frac{p}{N-p}\right)^p} \right]^{p'}=\left[1+\frac{1-\frac{\lambda}{N-1}}{1-\lambda\left(\frac{p}{N-p}\right)^p}-1\right]^{\frac{1}{\frac{1-\frac{\lambda}{N-1}}{1-\lambda\left(\frac{p}{N-p}\right)^p}-1}\left[\frac{1-\frac{\lambda}{N-1}}{1-\lambda\left(\frac{p}{N-p}\right)^p}-1\right]\frac{p}{p-1}}
$$
and, by l'Hospital's rule, 
$$
\lim_{p\to1^+}\left[\frac{1-\frac{\lambda}{N-1}}{1-\lambda\left(\frac{p}{N-p}\right)^p}-1\right]\frac{p}{p-1}=\frac{\lambda\left[\frac{N}{N-1}-\ln(N-1)\right]}{N-1-\lambda},
$$
we have that
$$
\lim_{p\to 1^+}\left[\frac{1-\frac{\lambda}{N-1}}{1-\lambda\left(\frac{p}{N-p}\right)^p} \right]^\frac{p}{p-1}=e^\frac{\lambda\left[\frac{N}{N-1}-\ln(N-1)\right]}{N-1-\lambda}.
$$
Consequently, by \eqref{EstUp0} and Young's inequality, $u_p$ is bounded in $BV(\Omega)$ and so, by compact embedding, up  to subsequences,
$$
u_p\to u\quad\mbox{ in } L^s(\Omega),\; s\in [1, 1^*)\quad\mbox{ as }\;p\to 1^+ .
$$
\fim

\medskip 

\subsection{Extreme case and beyond}\label{r61}
Here we collect some remark and examples that give some insight on the case in which  
$$
\|f\|_{L^N(\Omega)}S_N+\frac{\lambda}{N-1}\geq1.
$$

\begin{example}\label{r61b}\label{exa}
 Consider $f\geq 0$ not identically zero. Recall that, if $\|f\|_{L^N(\Omega)}S_N>1$ and $v_p\in W^{1,p}_0(\Omega)$ is a solution to the problem 
$$
\left\{
\begin{array}{rclr}
-\Delta_p v_p&=&f &\mbox{ in }\;\Omega,\\
v_p&=&0&\mbox{ on }\;\partial\Omega,
\end{array}
\right.
$$
in \cite{mst}, it is proved that
\begin{equation}
\label{DiVp}
\lim_{p\to 1^+}\int_{\Omega}|\nabla v_p|^p\,dx=\infty.
\end{equation}
On the other hand, taking $v_p$ as test function in \eqref{WSol}, we have
\begin{align*}
\int_{\Omega}|\nabla u_p|^{p-2}\nabla u_p\cdot\nabla v_p\,dx&=\lambda\int_{\Omega}\frac{1}{|x|^p} u_p^{p-1}v_p\,dx+\int_{\Omega}fv_p\,dx\\
&\geq \int_{\Omega} fv_p\,dx\\
&=\int_{\Omega}|\nabla v_p|^p\,dx.
\end{align*}
By H\"older inequality, this implies that
$$
\int_{\Omega}|\nabla v_p|^p\,dx\leq \int_{\Omega}|\nabla u_p|^p\,dx
$$
and so, by \eqref{DiVp}, 
$$
\lim_{p\to 1^+}\int_{\Omega}|\nabla u_p|^p\,dx=\infty.
$$
\end{example}
\begin{example}\label{exb} Let $u_p, v_p\in W^{1,p}_0(B_R(0))$ be solutions to the problems
$$
\label{EP1}
\left\{
\begin{array}{rclr}
-\Delta_p u &=& \frac{\lambda}{|x|^p}|u|^{p-2}u+1& \hbox{ in }\;B_R(0),\\
u&=&0&\hbox{ on }\; \partial B_R(0)
\end{array}
\right.
$$
and 
\begin{equation}
\label{EP2}
\left\{
\begin{array}{rclr}
-\Delta_p v &=& 1& \hbox{ in }\;B_R(0),\\
v&=&0&\hbox{ on }\; \partial B_R(0),
\end{array}
\right.
\end{equation}
respectively. Note that, if $R>N$,
$$
S_N\|1\|_{L^N(\Omega)}=\frac{1}{N |B_1(0)|^\frac{1}{N}} R|B_1(0)|^\frac{1}{N}=\frac{R}{N}>1.
$$
It is known that the explicit solution to Problem \eqref{EP2} (see \cite{Lindqvist1987, K}) is given by
$$
v_p(x)= v_p(0)\left(1-\left(\frac{|x|}{R}\right)^\frac{p}{p-1}\right),
$$
where
$$
v_p(0)=\left(\frac{1}{N}\right)^\frac{1}{p-1}\frac{p-1}{p}R^\frac{p}{p-1}.
$$

Note that, by weak maximum principle,
$$
u_p(x)\geq v_p(x)\qquad\mbox{ a.e. }\; x\in B_R(0).
$$
Since, for $R>N$ (see \cite{K}),
$$
\lim_{p\to 1^+}v_p(x)=+\infty\qquad\hbox{ a.e. }\; x\in B_R(0).
$$
Hence
$$
\lim_{p\to 1^+}u_p(x)=+\infty\qquad\hbox{ a.e. }\; x\in B_R(0).
$$

\end{example}

Now  let us  investigate the extreme case $\lambda=N-1$ and $f=0$.

Let $\lambda = N-1$ and  consider, for simplicity,  problem 
\begin{equation}
\label{Pf0}
\left\{
\begin{array}{rclr}
-\Delta_1 u &=& \frac{\lambda}{|x|}s(x)& \hbox{ in }\;B_1(0),\\
u&=&0&\hbox{ on }\; \partial B_1(0),
\end{array}
\right.
\end{equation}
where $s(x)\in \hbox{Sgn}(u(x))$ a.e. $x\in \Omega$. It is a matter of computation to show that any smooth decreasing function $v:(0,1]\mapsto \re^+$ satisfying $v(1)=0$ can give rise to a solution of problem \eqref{Pf0} in the sense of Definition \ref{DefSol}. In particular
$$
u(x)=\frac{1}{|x|^{\alpha}}-1
$$
is a solution of \eqref{Pf0} in the sense of Definition \ref{DefSol} provided $0<\alpha <N-1$. A part from reflecting once again the structural non-uniqueness feature of such an operator, this example also shows how solutions to our problem are not, in general, expected to be bounded.

\medskip 

\section{More general data and explicit examples }

\subsection{More general data}
In this subsection, we shall extended our previous results to the case of a general datum $f$, that is, when it belongs to the Lorentz space  $L^{N,\infty}(\Omega)$ (or the dual space $W^{-1,\infty}(\Omega)$). Let $f\in L^{N,\infty}(\Omega)$ not identically zero and $0<\lambda<N-1$, which satisfy the following condition
\begin{equation}
\label{H1}
\gamma\|f\|_{L^{N,\infty}(\Omega)}+\frac{\lambda}{N-1}\leq 1,
\end{equation}

where we recall that 
\begin{equation}
\label{ga}
\gamma = \left[(N-1)|B_1(0)|^\frac{1}{N}_N\right]^{-1}
\end{equation}
is the best constant in the embedding of $W^{1,1}(\Omega)\subset L^{\frac{N}{N-1},1}(\Omega)$ (see \eqref{bestl}). 
Note that, by same arguments to prove the existence of solution in \eqref{Pp} when $f\in L^N(\Omega)$, it can be prove that, if $f$ belongs to $L^{N,\infty}(\Omega)$, there exists $u_p\in W^{1,p}_0(\Omega)$ such that
$$
\int_{\Omega}|\nabla u_p|^{p-2}\nabla u_p\cdot \nabla  \varphi\,dx=\int_{\Omega}\frac{1}{|x|^p}|u_p|^{p-2}u_p\varphi\,dx+\int_{\Omega}f\varphi\,dx\qquad\varphi\in W^{1,p}_0(\Omega).
$$
Taking $u_p$ as test function above, using the Hardy and Sobolev inequalities, we obtain
$$
\int_{\Omega}|\nabla u_p|^p\,dx\leq
 \left[\frac{\gamma \|f\|_{L^{N, \infty}(\Omega)}}{1-\lambda\left(\frac{p}{N-p}\right)^p}\right]^\frac{p}{p-1}|\Omega|\,. 
$$
Thus, using Young inequality and  by \eqref{H1}, we can reason as before in order to get  the existence  $u \in BV(\Omega)$ such that, up least to subsequences,
\begin{eqnarray*}
u_p\to u&\mbox{ in }\;L^s(\Omega)\;s\in [1, 1^*),\\
u_p\rightharpoonup u&\mbox{ in }\; L^{1^*}(\Omega),\\
\nabla u_p\rightharpoonup u& *-\mbox{ weak in }\; BV(\Omega),
\end{eqnarray*}
as $p\to 1^+$. From now on, a line by line  extensions of the  arguments we used to prove Theorem \ref{exi}, it can be checked  that $u$ is a solution for \eqref{P} in sense of Definition \ref{DefSol}.

 \medskip 
 Analogously, if  $f \in {W^{-1,\infty}(\Omega)}$ is such that
$$
\label{H2}
\|f\|_{W^{-1,\infty}(\Omega)}+\frac{\lambda}{N-1}\leq 1,
$$
remembering that by Theorem \ref{exip} there are solutions $u_p$ of \eqref{Pp} that satisfy
 $$
\int_{\Omega}|\nabla u_p|^p\,dx\leq
  \left[\frac{ \|f\|_{W^{-1, \infty}(\Omega)}}{1-\lambda\left(\frac{p}{N-p}\right)^p}\right]^\frac{p}{p-1}, $$
  without any new effort one can deduce the existence of a distributional solution of   \eqref{P} in sense of Definition \ref{DefSol} as well.
   
   \begin{remark}
Let $f\in L^{N, \infty}(\Omega)$ be a nonnegative decreasing radially function. Consider the problem
$$
 \left\{
 \begin{array}{rclr}
 -\Delta_1 u&=&\frac{\lambda}{|x|}+f&\hbox{ in }\;\Omega\\
                 u&=&0&\hbox{ on }\;\partial\Omega.
  \end{array}
  \right.
$$
which, as in \cite{ct}, can be studied by means of the solutions of   
$$
 \left\{
 \begin{array}{rclr}
 -\Delta_p u_p&=&\frac{\lambda}{|x|}+f&\hbox{ in }\;\Omega\\
                 u_p&=&0&\hbox{ on }\;\partial\Omega,
  \end{array}
  \right.
$$
 as $p\to 1^+$. Moreover, under the condition \eqref{H1}, arguing as in \cite{ct}, we can conclude the same results as Theorem \ref{AsympB}.
   \end{remark}

\subsection{More explicit examples.} Recall that the existence of solutions to \eqref{P} is proved under the conditions 
\begin{equation}
\label{HE}
S_N\|f\|_{L^N(\Omega)}+\frac{\lambda}{N-1}\leq 1 \quad \hbox{ or }\quad \gamma\|f\|_{L^{N, \infty}(\Omega)}+\frac{\lambda}{N-1}\leq 1.
\end{equation}
We shall consider three  examples in case of generic $f$ and $\lambda$, inspired by the paper \cite{LatorreSegura2018, ct}, where these conditions are extremely satisfied or not. We start with a non-extreme case. 

\begin{example}  We assume {$$ \|f\|_{L^N(\Omega)}S_N+\frac{\lambda}{N-1}<1\,.$$}
Let $ a,\,\delta>0$ be such that
$$
\delta<\frac{1-\frac{\lambda}{N-1}}{|B_1(0)\setminus B_a(0)|^\frac{1}{N}S_N}.
$$
Let the function $f\>:\>[0,1]\to \mathbb{R}$ given by
$$
f(s)=\left\{
\begin{array}{lcr}
0& \quad&\mbox{if}\; s\in[0,a],\\
\delta&\quad&\mbox{if}\; s\in (a,1].
\end{array}
\right.
$$
Note that, $f(|x|)$, for all $x\in B_1(0)$ is well defined and it satisfies the strict inequality of the condition \eqref{HE}. We assume that the vector field ${\bf z}\in L^\infty(B_1(0),\mathbb{R}^N)$ is defined as
$$
{\bf z}(x)=\left\{
\begin{array}{lcr}
-\frac{\lambda}{N-1}\frac{x}{|x|}& \quad&\mbox{ if }\; x\in B_a (0)\setminus\{0\},\\
-\left(\frac{\lambda}{N-1}\frac{1}{|x|}+\frac{\delta}{N}\right)x&\quad&\mbox{if}\; x\in B_1\setminus B_a.\\
\end{array}
\right.
$$
and $1= s(x)\in \mbox{Sgn}(u(x))$ a.e. on $B_1(0)$.  Consequently $u=0$ in $B_1(0)$ and $u=0$ on $\partial B_1(0)$ is a solution   for \eqref{P} in the sense of the Definition \ref{DefSol}. 
\end{example}

In the following example we show the existence of a non-trivial solution in the extremal case.  
\begin{example}{We consider $$\gamma\|f\|_{L^{N, \infty}(\Omega)}+\frac{\lambda}{N-1}=1\,. $$}
Let $0<\alpha< 1$. Consider the problem 
$$
\label{NT}
\left\{  
\begin{array}{rclr}
-\Delta_pu &=&\frac{\lambda}{|x|}+\frac{(1-\alpha)(N-1)}{|x|}&\hbox{ in }\; B_1(0)\\
        	      u &=&0&\hbox{ on }\; \partial B_1(0),
\end{array}
\right.
$$
where $\lambda=\alpha(N-1)$. Note that, for $\gamma=\left[(N-1)|B_1(0)|^\frac{1}{N}_N\right]^{-1}$,
$$
\label{EH}
\gamma \left\|\frac{(1-\alpha)(N-1)}{|x|}\right\|_{L^{N, \infty}(\Omega)}+\frac{\lambda}{N-1}=1.
$$
 The solution in $W^{1,p}_0(\Omega)$ to the problem above is given by
 $$
u_p(x)=1-|x|, 
$$
which clearly converges to $1-|x|$ a.e. $x\in\Omega$ as $p\to 1^+$. But this function together with 
$$
{\bf z}=-\frac{x}{|x|}\quad\hbox{ a.e}\;x\in \Omega\quad\hbox{and }\quad s(x)=1\quad\hbox{ a.e.} \;x\in\Omega
$$
satisfy the conditions (1)-(3) of Definition \ref{DefSol}.
\end{example}

  \medskip
  
 We conclude with the following example, in which one also get that  the smallness assumption \eqref{H1}, for fixed $0<\lambda<N-1$,  is  sharp.
 \begin{example}
Let $f=(a-\lambda)/|x|$ with $\lambda<a\leq N-1$ and $\Omega=B_1(0)$. We consider the following approximate $p-$Laplacian problem  
\begin{equation}
\label{PLp}
\left\{
\begin{array}{rclr}
-\Delta_p u&=&\frac{\lambda}{|x|}+f&\mbox{ in }\; B_1(0),\\
u&=&0&\mbox{on }\;\partial B_1(0).
\end{array}
\right.
\end{equation}
A solution to the problem \eqref{PLp} is given by
$$
u_p(x)=(1-|x|)\left(\frac{a}{N-1}\right)^\frac{1}{p-1}\qquad\mbox{ for all }\; x\in B_1(0).
$$
Since
\begin{equation}
\label{Sol}
\lim_{p\to1^+}u_p(x)=
\left\{
\begin{array}{lcl}
0&\hbox{ if }& a<N-1\\
1-|x|&\hbox{ if }&a=N-1\\
+\infty&\hbox{ if }&a>N-1,
\end{array}
\right.
\end{equation}
we may infer that $0$ and $1-|x|$ are solutions (in sense of Definition \ref{DefSol}) to the problem
\begin{equation}
\label{PL}
\left\{
\begin{array}{rclr}
-\hbox{div}\,{\bf z}&=&\frac{\lambda}{|x|}s(x)+f&\mbox{ in }\; B_1(0),\\
u&=&0&\mbox{on }\;\partial B_1(0),
\end{array}
\right.
\end{equation}
where ${\bf z}(x)=-\frac{x}{|x|}$ and $s(x)=1\in \hbox{Sgn}(u(x))$ a.e. $x\in \Omega$. Furthermore, \eqref{Sol} implies that \eqref{H1} is sharp.
\end{example}

{\bf Acknowledgment.}
Francesco Petitta is partially supported by the Gruppo Nazionale per l’Analisi Matematica, la Probabilità e le loro Applicazioni (GNAMPA) of the Istituto Nazionale di Alta Matematica (INdAM). 
Juan C. Ortiz Chata is supported by FAPESP 2021/08272-6 and 2022/06050-9, Brazil.

\section*{Conflict of interest declaration}

The authors declare no competing interests.


\begin{thebibliography}{99}



\bibitem{Alvino1977}{\sc Alvino, A.:} {\em Sulla disegualianza di {S}obolev in spazi di {L}orentz}, Bolletino dell'Unione Matematica Italiana, 14-A(5) (1977) no. 1, 148--156.

\bibitem{AmbrosioFuscoPallara}{\sc Ambrosio, L., Fusco, N., Pallara, D.:} {\em Functions of bounded variation and free discontinuity problems}, Oxford University Press, Oxford (2000).

\bibitem{acm2001}{\sc Andreu, F.; Ballester, C.; Caselles, V.; Maz\'on, J.M.:} {\em The {D}irechlet problem for the total variation flow}, J. Funct. Anal. 180 (2001), no. 2, 347--403.


\bibitem{AndreuMazonMollCaselles2004}{\sc Andreu, F.; Maz\'on, J.M.; Moll, S.; Caselles, V.:} {\em The minimizing total variation flow with measure initial conditions}, Comm. Contemp. Math. 6 (2004), no. 3, 431--494.

\bibitem{AndreuCasellesMazon2004}{\sc Andreu, F., Caselles, V. and Maz\'on, J.M.:} {\em Parabolic quasilinear equations minimizing linear growth functionals}, {P}rogress in mathematics 223, {B}irkh\"auser  {V}elarg, {B}asel, 2004.

\bibitem{Anzellotti1983} {\sc Anzellotti, G.:} {\em Pairings between measures and bounded functions and compensated compactness}, Ann. Mat. Pura Appl., 135 (1983),  no. 1, 293--318.


\bibitem{abo} {\sc Arcoya, D.; Boccardo, L.; Orsina, L.:} {\em Hardy potential versus lower order terms in Dirichlet problems: regularizing effects}, Math. Eng. 5 (2023), no. 1, Paper No. 004, 14 pp.

\bibitem{AttouchButtazzoMichaille2006} {\sc Attouch, H.; Buttazzo, G.; Michaille, G.:} {\em Variational analysis in Sobolev spaces: applications to PDEs and optimization.} MPS-SIAM, Philadelphia, (2006). 

\bibitem{BCRS} {\sc  Bertalmio M.; Caselles V.; Roug\'e B.;   Sol\'e A. }: 
TV based image restoration with local constraints.
Special issue in honor of the sixtieth birthday of Stanley Osher, 
J. Sci. Comput. 19 (2003), No. 1-3, 95--122.

\bibitem{boc} {\sc Boccardo, L.:}, {\em Hardy potentials and quasi-linear elliptic problems having natural growth terms},
Progr. Nonlinear Differential Equations Appl.,
Birkhäuser Verlag, Basel, 63 (2005), 67--82.

\bibitem{BoccardoMurat1992}{\sc Boccardo, L.; Murat, F.:} {\em Almost everywhere convergence of the gradients of solutions to elliptic and parabolic equations}, Nonlinear Anal. 19 (1992), no. 19, 581--597.

\bibitem{BoccardoOrsinaPeral2006}{\sc Boccardo, L; Orsina, L; Peral, I.:} { \em A remark on existence and optimal summability of solutions of elliptic problems involving Hardy potential}, Discrete Contin. Dyn Syst. 16 (2006), no. 3, 513--523.

\bibitem{CKN1984}{\sc Caffarelli, L.; Kohn, R.; Nirenberg, L.:} {\em First order interpolation inequalities with weights}, Compositio Math. 53 (1984), no. 3, 259--275.

\bibitem{Caselles2011}{\sc Caselles, V.:} {\em On the entropy conditions for some flux li equations}, J. Differential Equations, 250 (2011), no. 8, 3311--3348.

\bibitem{ChenFried1999}{\sc Chen, G.-Q; Frid, H.:} {\em Divergence-measure fields and hyperbolic conservation laws}, Arch. Ration. Mech. Anal. 147 (1999), no. 2, 89--118.

\bibitem{ct}{\sc Cicalese, M.; Trombetti, C.:} {\em Asymptotic behaviour of solutions to {$p$}-{L}aplacian equation}, Asympt. Anal. 35 (2003), no. 1, 27--40.



\bibitem{DAnconaFanelli2007}{\sc D'Ancona, P.; Fanelli, L.:} {\em Decay Estimates for the Wave and Dirac Equations with a Magnetic Potential}, Comm. Pure Appl. Math, 60 (2007), no. 3, 357--392.



{ \bibitem{dgop} {\sc De Cicco V.; Giachetti, D.; Oliva, F.; Petitta, F.: }{\em The Dirichlet problem for singular elliptic
equations with general nonlinearities}, Calc. Var. Partial Differential Equations, 58 (2019), no. 4, 1--40.\bk }

\bibitem{Demengel2004}{\sc Demengel, F.:} {\em Functions locally almost 1-harmonic}, App. Anal. 83 (2004), no. 9, 865--896.

\bibitem{Demengel2005}{\sc Demengel, F.:} {\em Some existence results for noncoercive "1-Laplacian" operator}, Asymptot. Anal. 43 (2005), no. 4, 287--322.

\bibitem{EvansGariepy}{\sc Evans, L.C.; Gariepy, R.F.:} {\em Measure theory and fine properties of functions}, Textbook in Matematics,  CRC Press, Boca Raton, FL (2015), xiv+299.

\bibitem{AzoreroAlonso1998}{\sc Garc\'{\i}a Azorero, J. P.; Peral Alonso, I.:} {\em Hardy inequalities and some critical elliptic and parabolic problems}. J. Differential Equations, 144 (1998), no. 2, 441--476.

\bibitem{LMP2011}{ \sc Leonori, T.; Mart\'{\i}nez-Aparicio, P.J.; Primo, A.:} {\em Nonlinear elliptic equations with {H}ardy potential and lower order term with natural growth}, Nonlinear Anal. 74 (2011), no. 11, 3556--3569.

\bibitem{lp} {\sc Leonori T.,  Petitta F. :} {\em Existence and regularity results for some singular elliptic problems}, Adv. Nonlinear Stud., 7 (2007), n. 3, 329--344


\bibitem{Lindqvist1987}{\sc Lindqvist, P.:} {\em Stability for the solutions of {${\rm div}\,(|\nabla u|^{p-2}\nabla u)=f$} with varying {$p$}}, J. Math. Anal. Appl. 127 (1987), no. 1, 93--102.

    \bibitem{K}  {\sc  Kawohl B.}: {\em On a family of torsional creep problems.}
     J. Reine Angew. Math.  410 (1990),  1--22.
     
     \bibitem{ka} {\sc Kawohl B.}:  {\em From $p$-Laplace to mean curvature operator and related questions.} In : Progress in Partial Differential Equations: the Metz Surveys, Pitman Res. Notes Math. Ser.,  249 (1991), pp. 40--56.
     

\bibitem{Kawohl2007}{\sc Kawol, B; Schuricht, F.:} {\em Dirichlet problems for the 1-{L}aplace operator, including the eigenvalue problem}. Commun. Contemp. Math. 9 (2007), no. 4, 515--543.

\bibitem{LatorreSegura2018}{\sc Latorre, M.; Segura de Le\'{o}n, S.:} {\em Existence and comparison results for an elliptic equation involving the $1-${L}aplacian and {$L^1$}-data}. J. Evol. Equ., 18 (2018), no. 1, 1--28.

\bibitem{lops}{\sc Latorre, M.; Oliva F.; Petitta F. ;Segura de Le\'{o}n, S.:} {\em The Dirichlet problem for the 1-Laplacian with a general singular term and  L1 -data}
Nonlinearity, 34 (2021), no. 3, 1791--1816.

{ \bibitem{mst}{\sc Mercaldo, A.; Segura de Le\'on, S.; Trombetti, C.:} {\em On the behaviour of the solutions to $p$-Laplacian equations as $p$ goes to 1}. Publ. Mat. 52 (2008), no. 2, 377--411.\bk }

\bibitem{MST2009}{\sc Mercaldo, A.; Segura de Le\'{o}n, S.; Trombetti, C.:}{ \em On the solutions to 1-{L}aplacian equation with {$L^1$} data}, J. Funct. Anal. 256 (2009), no. 8, 2387--2416.

\bibitem{M} {\sc  Meyer Y.:}  {Oscillating Patterns in Image Processing and Nonlinear Evolution Equations: The Fifteenth Dean Jacqueline B. Lewis Memorial Lectures}.   Providence, RI: American Mathematical Society, 2001.



\bibitem{ONeil1968}{\sc O'Neil, R.:} {\em Integral transforms and tensor products on {O}rlicz spaces and $L(p, q)$ spaces}, J. Anal. Math. 21 (1968), 1--276.

\bibitem{OrsinaOrtiz2024}{\sc Orsina, L.; Ortiz Chata, J.C.:} {\em Hardy potential in existence and optimal summability of solutions of nonlinear elliptic problem}. Preprint.

\bibitem{OrtizPimentaSegura2021}{\sc Ortiz Chata, J. C.; Pimenta, M.T.O.; Segura de Le\'on, S.:} {\em Anisotropic 1-{L}aplacian problems with unbounded weights}, NoDEA Nonlinear Differential Equations Appl., 28 (2021), no. 57, 40.

\bibitem{OsSe}
{\sc Osher S.; Sethian J.:} Fronts propagating with curvature-dependent speed: algorithms based on Hamilton-Jacobi formulations. Journal of Computational Physics, 79 (1998), no.1, 12--49. 

 \bibitem{ROF}{\sc   Rudin L.I. Osher S. and Fatemi  E.}: {\it Nonlinear total variation based noise removal algorithms}, 
   Physica D. {\bf 60} (1992), 259--268.

    \bibitem{Sapiro}  {\sc  Sapiro G.:}  {Geometric partial differential equations and image analysis,} Cambridge University Press, 2001.

\bibitem{Talenti}{\sc Talenti, G.:} {\em Best constant in {S}obolev inequality}, Ann. Mat. Pura Appl. 110 (1976), 353--372.

\bibitem{Ziemer1989}{\sc Ziemer, W.P.:} {\em Weakly differentiable functions}, Spring-Velarg, New York, 120 (1989), xvi+308.







\end{thebibliography}
\end{document}